ORIGINAL PAPER

# On the Stability of IMEX Upwind gSBP Schemes for 1D Linear Advection-Diffusion Equations

**Sigrun Ortleb**[1]

## Abstract

A fully discrete energy stability analysis is carried out for linear advection-diffusion problems discretized by generalized upwind summation-by-parts (upwind gSBP) schemes in space and implicit-explicit Runge-Kutta (IMEX-RK) schemes in time. Hereby, advection terms are discretized explicitly, while diffusion terms are solved implicitly. In this context, specific combinations of space and time discretizations enjoy enhanced stability properties. In fact, if the first- and second-derivative upwind gSBP operators fulfill a compatibility condition, the allowable time step size is independent of grid refinement, although the advective terms are discretized explicitly. In one space dimension it is shown that upwind gSBP schemes represent a general framework including standard discontinuous Galerkin (DG) schemes on a global level. While previous work for DG schemes has demonstrated that the combination of upwind advection fluxes and the central-type first Bassi-Rebay (BR1) scheme for diffusion does not allow for grid-independent stable time steps, the current work shows that central advection fluxes are compatible with BR1 regarding enhanced stability of IMEX time stepping. Furthermore, unlike previous discrete energy stability investigations for DG schemes, the present analysis is based on the discrete energy provided by the corresponding SBP norm matrix and yields time step restrictions independent of the discretization order in space, since no finite-element-type inverse constants are involved. Numerical experiments are provided confirming these theoretical findings.

✉ Sigrun Ortleb
ortleb@mathematik.uni-kassel.de

[1] Department of Mathematics and Natural Sciences, University of Kassel, Untere Königsstraße 86, 34109 Kassel, Germany

## 1 Introduction

Advection-diffusion problems describe a wide range of phenomena arising due to the combination of two different physical processes, i.e., advective transport of a quantity due to a given velocity field and movement of particles from areas of higher to lower concentration. For instance, advection-diffusion equations arise as a model for the time evolution of chemical or biological species in a flowing medium. If the resulting partial differential equations are discretized on refined spatial grids, the presence of diffusion fluxes usually generates increasingly stiff problems. Particularly, explicit time discretization of semi-discrete advection-diffusion problems often results in an increasingly severe grid-dependent time step restriction since the time step $\Delta t$ scales with the grid length scale $\Delta x$ as $\Delta t = \mathcal{O}(\Delta x^2)$. However, purely implicit time discretization requires the solution of large nonlinear systems of equations. Therefore, hybrid time integration schemes such as implicit-explicit (IMEX) methods have frequently been considered. Hereby, a common approach is to discretize the advective terms explicitly, while the diffusion terms are treated implicitly, see, e.g., [4, 5, 9, 37]. The additional complexity of a partially implicit treatment is justified by the observation that even though many applications contain nonlinear advection terms, the diffusion terms are often linear, only demanding the iterative solution of linear algebraic systems which are positive definite, symmetric, and sparse. In addition, the nonlinear systems resulting from potentially nonlinear diffusion terms are generally more efficient to solve than those arising from nonlinear advection terms.

While this approach alleviates the extremely restrictive time step scaling for the diffusion terms, it does not a priori dispose of the classical CFL-type time step restriction of the form $\Delta t = \mathcal{O}(\Delta x)$ due to the explicit time discretization of the advection terms. However, for specific combinations of space and time discretizations, enhanced stability properties may be achieved. For instance, Calvo et al. [9] designed specific IMEX Runge-Kutta (IMEX-RK) methods which guarantee grid-independent time step restrictions of the form $\Delta t = \mathcal{O}(c/a^2)$ for Fourier spatial discretizations, where $a$ and $c$ denote the advection and diffusion coefficients, respectively. Wang et al. [38] proved a similar unconditional $L^2$-stability result in the context of spatial discretization by the discontinuous Galerkin (DG) [12] method applied to linear advection-diffusion equations, with an extension to nonlinear problems in [39] and to the multi-dimensional case in [40]. Hereby, diffusion terms were discretized by the local discontinuous Galerkin (LDG) scheme [11]. In [27], an extension was provided for DG methods in one space dimension incorporating a larger class of diffusion schemes, namely, the $(\sigma, \mu)$ family [20, 21] of diffusion discretizations which also includes the Bassi-Rebay schemes [6, 7]. Furthermore, Fu and Shu [14] proved the grid-independent $L^2$-stability for the discretization of the diffusion terms by the embedded DG method and in [41], the direct DG method was considered. A related approach which provides the unconditional stability even if part of the problem is discretized explicitly was taken in [32, 33] designing unconditionally stable IMEX linear multistep schemes, where both the IMEX splitting and the employed multistep schemes fulfill specific properties.

In this work, upwind generalized summation-by-parts (gSBP) space discretizations for advection-diffusion problems will be considered with respect to the combination with IMEX time integration and regarding stable time step choices. These upwind gSBP schemes do not necessarily provide a replacement for the DG space discretization, but rather a broader framework.

This general framework of SBP schemes provides structure-preserving spatial discretizations also including DG type methods. SBP operators fulfill a mimetic property, in particular mimicking the analytical concept of integration-by-parts. Originally, discrete derivative operators with an SBP property have been introduced for finite difference schemes for various types of PDEs including advection-diffusion equations and the linearized compressible Navier-Stokes equations, see [10, 19, 25, 31, 35, 36]. Extensions of the SBP methodology have been provided for tensor-product grids on curvilinear elements [29] as well as for the general multidimensional case [18] including simplex elements. Via $L^2$ energy estimates, the spatial SBP operators automatically yield stable semi-discrete schemes for periodic solutions of a broad class of linear equations. In addition, they have substantially profited from a combination with weakly enforced boundary conditions, most prominently using simultaneous approximation terms (SATs) which were first developed in [10]. More recent investigations also focus on SBP operators within various popular classes of numerical schemes, e.g., finite volume schemes on unstructured dual grids [24], DG schemes with Legendre-Gauss-Lobatto nodes [15], DG schemes on Legendre-Gauss nodes both in one space dimension and on tensor-product grids [26], and flux-reconstruction schemes [28]. In fact, the definition in [30] allows to derive gSBP properties for nodal DG schemes in one space dimension on arbitrary nodal sets, as long as their corresponding quadrature rule is sufficiently exact to mimic integration-by-parts at the discrete level. An extension of the SBP framework departing from their central character was originally proposed by Mattsson [23], where upwind SBP schemes of finite difference type were constructed. These schemes consist of dual-pair SBP operators with non-central difference stencils. The incorporated upwind differencing leads to a built-in artificial dissipation in the semi-discrete energy analysis of periodic problems and has been used to construct upwind SBP schemes based on flux splitting for the shallow water equations in [22] and for general scalar hyperbolic conservation laws in [34].

Connections between nodal DG schemes and SBP methods have mainly fostered the discretization of conservation laws in the skew-symmetric form to preserve secondary quantities and to move towards semi-discrete energy or entropy stability results. In this context, the SBP properties of *cellwise* DG discretizations at the local level of one element represent the key aspect. In this work, emphasis is put on the *global* formulation of DG schemes in one space dimension which may be regarded as upwind gSBP schemes for a family of numerical flux functions ranging from upwind to central ones. In addition, we will take advantage of the connection to upwind gSBP schemes to improve the fully discrete stability estimates pertaining to DG and upwind gSBP schemes combined with the IMEX time integration.

While discretizing advection terms with a DG scheme is straightforward, various approaches to discretize diffusion terms have been introduced in the literature, since the discretization of higher order spatial derivatives is less natural. Either specially designed penalty terms are introduced as in [2, 8] or the advection-diffusion equation is rewritten into a system of first-order equations using auxillary variables for the solution derivatives as in [3, 6, 7, 11]. The first method by Bassi and Rebay [6], termed the BR1 scheme, has actually been the first extension of the DG scheme to the compressible Navier-Stokes equations. It is based on rewriting the viscous terms into a larger, extended first-order degenerate system of PDEs with the gradient as a new unknown. After this reformulation, the standard DG approach is applied to the extended system which necessitates to prescribe two types of numerical fluxes. Hereby, the BR1 scheme is the simplest approach, using arithmetic means for both types of fluxes. Cockburn and Shu [11] analyzed various methods based on the reformulation into a first-order PDE and derived conditions on the numerical fluxes to guarantee the stability, the convergence, and a suboptimal error estimate of order $N$ when using an approximation space of polynomial degree $N$. The analysis shows the sub-optimal convergence of the BR1 scheme for odd $N$, while the choice of alternating numerical fluxes usually associated with the LDG scheme by Cockburn and Shu [11] leads to the optimal convergence of order $N+1$. Being simple to code, parameter-free, and generic for nonlinear viscous fluxes and arbitrary grids, the BR1 scheme is still very popular and satisfies certain structural properties. In fact, in [17], the neutral behaviour of BR1 with respect to artificial dissipation over element interfaces and the resulting stability for the compressible Navier-Stokes equations has been proven.

The energy stability analysis in [27] for linear advection-diffusion problems discretized by the DG scheme in space using a $(\sigma, \mu)$ diffusion discretization has shown that the simple BR1 approach coupled with upwind fluxes for advection does not inherit the enhanced stability properties in the case of IMEX time integration. In this regard, this current work provides a remedy, since the DG scheme based on BR1 fluxes may be regarded as a specific second-derivative upwind gSBP scheme compatible with central fluxes for advection. Furthermore, unlike previous discrete energy stability investigations for DG schemes, the present analysis is based on the discrete energy provided by the corresponding SBP norm matrix and yields time step restrictions independent of the discretization order in space, since no FE type inverse constants are involved.

This paper is organized as follows. Section 2 introduces the properties of first- and second-derivative upwind gSBP operators and reviews their connection to DG schemes on the global level. Furthermore, some properties and estimates required for the fully discrete stability analysis are provided. In Sect. 3, an energy stability analysis based on the discrete energy defined by the associated upwind gSBP norm matrix is carried out for first- and second-order IMEX-RK time integrators combined with upwind gSBP schemes of arbitrary order in space. Numerical results verifying the conducted theoretical analysis are presented in Sect. 4. Finally, Sect. 5 contains a conclusion and an outlook on future work.

## 2 Discontinuous Space Discretization of the Linear advection-diffusion Equation

Throughout this work, we consider the linear advection-diffusion equation:

$$u_t + au_x = cu_{xx}, \qquad (x,t) \in (x_a, x_b) \times (0,T) \tag{1}$$

with the advective velocity $a > 0$ and the diffusion coefficient $c > 0$, supplemented by the periodic initial condition $u(x,0) = u_0(x)$ in $L^2(\Omega)$ and periodic boundary conditions.

For space discretization, the spatial domain $\Omega = (x_a, x_b)$ is partitioned into cells $I_i = (x_i, x_{i+1})$, $i = 1, \cdots, K$ with $x_1 = x_a$, $x_{K+1} = x_b$. We consider quasi-uniform grids with cell lengths $\Delta x_i = x_{i+1} - x_i$ satisfying $\frac{\Delta x_i}{\Delta x} \geqslant \rho \in \mathbb{R}^+$ for $\Delta x = \max_i \Delta x_i$ under grid refinement $\Delta x \to 0$. On the closure $\bar{I}$ of the reference interval $I = (-1, 1)$, a nodal set of solution points is specified and transferred to each grid cell, not necessarily including the cell boundaries in the nodal set. Furthermore, the discontinuities space discretization allows for a discontinuous representation of the semi-discrete solution at cell interfaces.

### 2.1 Discontinuous Galerkin (DG) Formulation

DG schemes allow for piece-wise smooth approximate solutions with potential discontinuities across element interfaces. Due to its flexibility and generality, the DG scheme is a popular numerical method in a variety of applications. The main advantages are its local conservation property, an arbitrarily high order of accuracy, and superconvergence capabilities. Dispensing with global continuity of the approximate solution, the DG approach results in relatively compact stencils which greatly facilitates both *hp*-adaptivity of the method and its implementation in parallel hardware environment. The DG scheme is built upon a variational formulation of the advection-diffusion equation. The basis functions and test functions used to define the DG scheme in one space dimension are taken from the finite element space $V_h = \{v \in L^2(\Omega) \mid v|_{I_i} \in \mathcal{P}^N(I_i),\ \forall i = 1, \cdots, K\}$, for the spatial domain $\Omega = (x_a, x_b) = \bigcup_{i=1}^{K} I_i$, where $\mathcal{P}^N(I_i)$ denotes the space of polynomial functions on $I_i$ of degree at most $N$. As usual in DG schemes, the functions in $V_h$ may be discontinuous across element boundaries.

For the discretization of diffusion terms, the PDE (1) will be reformulated as a first-order system of equations. Thus, the structure of the corresponding DG scheme is similar to the case of purely hyperbolic problems. The DG formulation for the linear diffusion equation will then be given in Sect. 2.2.2. To formulate the DG scheme in matrix-vector notation similar to SBP schemes, it is sufficient to consider a scalar hyperbolic conservation law in one space dimension given by

$$\frac{\partial}{\partial t}u(x,t) + \frac{\partial}{\partial x}f(u(x,t)) = 0, \quad (x,t) \in (x_a, x_b) \times (0, T). \tag{2}$$

The DG scheme constructs an approximation $u_h$ of $u$ which is piecewise continuous of the form

$$u_h(x,t)|_{I_i} = u_h^i(x,t) = \sum_{j=1}^{N+1} u_j^i(t) L_j^i(x) \tag{3}$$

using polynomial basis functions $L_j^i \in \mathcal{P}^N\left(I_i\right)$. Here, nodal DG schemes are considered, where the basis functions are the Lagrange polynomials corresponding to a set of solution points $\xi_\nu \in \bar{I},\ \nu = 1, \cdots, N+1$ which are the nodes $\xi_k$ of a suitable quadrature rule on the closure of the reference cell $(-1, 1)$. We assume that the quadrature rule is symmetric and exact for polynomials of degree $2N-1$ and denote the corresponding weights by $\omega_\nu,\ \nu = 1, \cdots, N+1$. More precisely, we consider the Lagrange polynomials defined by $L_j^i(\Lambda_i(\xi_\nu)) = \delta_{\nu j}$, where $\delta_{\nu j}$ denotes the Kronecker delta and $\Lambda_i$ transforms the reference cell $I$ to the specific sub-interval $I_i$, i.e., $\Lambda_i(\xi) = \xi\frac{x_{i+1}-x_i}{2} + \frac{x_i+x_{i+1}}{2}$.

The DG scheme in a weak form, obtained by multiplying (2) by a test function $L_k^i$, integrating in space over a cell $I_i$, and applying partial integration, is given by

$$\frac{\mathrm{d}}{\mathrm{d}t}\int_{I_i} u_h L_k^i \mathrm{d}x + f_i^* L_k^i(x_i) - f_{i+1}^* L_k^i(x_{i+1}) - \int_{I_i} f_h \frac{\mathrm{d}L_k^i}{\mathrm{d}x}\mathrm{d}x = 0 \tag{4}$$

with $f_i^* = f^*\left(u_h^{i-1}(x_i,t), u_h^i(x_i,t)\right)$ denoting the values of a consistent numerical flux function $f^*$ evaluated at the left- and right-hand side limits of $u_h$ at the cell boundaries $x_i$ and the flux polynomial $f_h(x,t)|_{I_i} = \sum_{j=1}^{N+1} f(u_j^i) L_j^i(x)$, based on the pointwise values $u_j^i = u_h\left(\Lambda_i(\xi_j), t\right)$.

The integrals in (4), are then solved numerically by the chosen quadrature rule which exactly evaluates the last integral due to the degree of exactness.

The DG scheme in a strong from is obtained by a second partial integration of (4) resulting in the variational formulation

$$\frac{\mathrm{d}}{\mathrm{d}t}\int_{I_i} u_h L_k^i \mathrm{d}x + \int_{I_i} \frac{\partial f_h}{\partial x} L_k^i \mathrm{d}x = \left[f_i^* - f_h^i(x_i)\right] L_k^i(x_i) - \left[f_{i+1}^* - f_h^i(x_{i+1})\right] L_k^i(x_{i+1}). \tag{5}$$

Using the matrix-vector notation as shown, e.g., in [1, 26], (5) rewrites in a simpler form when transformed to the reference cell $I$. For this purpose, we define the matrices $\widehat{\mathbf{M}}$ and $\widehat{\mathbf{S}}$ by

$$\widehat{M}_{jk} = \sum_l \omega_l L_j(\xi_l) L_k(\xi_l) = \widehat{M}_{kj} \approx \int_{-1}^{1} L_j L_k \mathrm{d}\xi,$$

$$\widehat{S}_{jk} = \sum_l \omega_l L_j(\xi_l) L_k'(\xi_l) = \int_{-1}^{1} L_j L_k' \mathrm{d}\xi,$$

where $L_k : \bar{I} \to \mathbb{R}$ are the Lagrange polynomials corresponding to the quadrature nodes $\xi_k$. From the Lagrange interpolation property $L_j(\xi_l) = \delta_{jl}$, we directly obtain

$$\widehat{M}_{jk} = \omega_j \delta_{jk}, \; \widehat{S}_{jk} = \omega_j L_k'(\xi_j), \tag{6}$$

hence $\widehat{\mathbf{M}}$ is diagonal. Using the solution vector $\mathbf{u}^i$ and the vector of flux values $\mathbf{f}^i$,

$$\mathbf{u}^i = (u_1^i, \cdots, u_{N+1}^i)^{\mathrm{T}} \text{ with } u_j^i = u_h\big(\Lambda_i(\xi_j), t\big),$$

$$\mathbf{f}^i = (f_1^i, \cdots, f_{N+1}^i)^{\mathrm{T}} \text{ with } f_j^i = f(u_j^i),$$

as well as the transformation of $f_h(x,t)|_{I_i}$ to the reference cell by $f_h^i(\xi,t) = \sum_{j=1}^{N+1} f(u_j^i) L_j(\xi)$, the abbreviations $f^{*,i}(1) = f_{i+1}^*$ and $f^{*,i}(-1) = f_i^*$ and the collection of the basis functions in the vector valued function $\mathbf{L}(\xi) = (L_1(\xi), \cdots, L_{N+1}(\xi))^{\mathrm{T}}$, (5) is equivalent to the matrix-vector formulation

$$\frac{\Delta x_i}{2} \widehat{\mathbf{M}} \frac{\mathrm{d}\mathbf{u}^i}{\mathrm{d}t} + \widehat{\mathbf{S}}\mathbf{f}^i = \left[\big(f_h^i - f^{*,i}\big)\mathbf{L}\right]_{-1}^{1},$$

see [1, 26]. By definition, the matrix $\widehat{\mathbf{M}}$ is invertible. Multiplying by the inverse of $\widehat{\mathbf{M}}$, we obtain the final matrix-vector formulation of the DG scheme

$$\frac{\Delta x_i}{2} \frac{\mathrm{d}\mathbf{u}^i}{\mathrm{d}t} + \widehat{\mathbf{D}}\mathbf{f}^i = \widehat{\mathbf{M}}^{-1}\left[\big(f_h^i - f^{*,i}\big)\mathbf{L}\right]_{-1}^{1}, \tag{7}$$

where the entries of $\widehat{\mathbf{D}} = \widehat{\mathbf{M}}^{-1}\widehat{\mathbf{S}}$ can be derived from (6) and are given by $\widehat{D}_{jk} = L_k'(\xi_j)$.

**Remark 1** It is worth noting that the matrix $\widehat{\mathbf{D}}$ provides a local first-derivative gSBP operator on a specific DG cell, see, e.g., [26], but outsources the interface fluxes to be dealt with similar to the weak boundary treatment of SATs. In contrast, the next section investigates the upwind gSBP structure of the global DG scheme including the presence of numerical fluxes on cell interfaces.

## 2.2 Upwind gSBP Space Discretization

Upwind SBP schemes, first derived in [23], are a class of finite difference methods which combine the structure-preserving properties of SBP derivative operators with artificial dissipation by choosing upwind directions. These schemes consist of dual-pair SBP operators with non-central difference stencils which lead to a built-in artificial dissipation when combined with flux splitting. While traditionally, SBP schemes were set up on equidistant grid points, we will consider generalized versions in this work on non-uniform nodal sets not necessarily containing boundary points which have been termed gSBP schemes. These schemes may be further modified by including upwind stencils, resulting in upwind gSBP schemes which will be considered in this work. The common goal of schemes within the SBP framework is to transfer continuous energy estimates to the semi-discrete level. For instance, considering solutions $u(x, t)$ of the linear advection equation:

$$u_t + au_x = 0, \ a > 0, \quad (x, t) \in (\alpha, \beta) \times (0, T) \tag{8}$$

subject to the homogeneous boundary condition $u(x_a, t) = 0$, the continuous energy estimate

$$\frac{\mathrm{d}}{\mathrm{d}t} \|u(x, t)\|^2_{L^2(\alpha,\beta)} + au^2(\beta, t) = 0 \tag{9}$$

can be proven. By the SBP property, this estimate may be transferred to the spatial discretization. The corresponding operators approximate the first derivative $\frac{\partial}{\partial x}$ in (8) and are based on a set of grid points $\{x_i\}_{1 \leqslant i \leqslant M} \subset [\alpha, \beta]$ called a nodal set. On this nodal set, approximate solutions of (8) are given by time-dependent solution vectors

$$\mathbf{u}(t) = (u_1(t), \cdots, u_M(t))^{\mathrm{T}} \approx (u(x_1, t), \cdots, u(x_M, t))^{\mathrm{T}},$$

and $u_x$ in (8) is replaced by the term $\mathbf{D}^- \mathbf{u}$, utilizing an upwind gSBP operator $\mathbf{D}^-$ as defined below.

**Definition 1** A pair of difference operators denoted by $\mathbf{D}^+$ and $\mathbf{D}^-$, both approximating the first derivative $\frac{\partial}{\partial x}$ on an interval $(\alpha, \beta)$, are called upwind gSBP operators of degree $q$ if

i. the matrices $\mathbf{D}^+$ and $\mathbf{D}^-$ provide accurate approximations to $\frac{\partial}{\partial x}$ of degree $q$, i.e.,

$$\mathbf{D}^{\pm} \mathbf{x}^k = k\mathbf{x}^{k-1}, \quad 0 \leqslant k \leqslant q,$$

where $\mathbf{x}^k = \left(x_1^k, \cdots, x_M^k\right)^{\mathrm{T}}$ is the representation of the monomials $x^k$ on the grid points;

ii. there exists a positive definite symmetric matrix $\mathbf{M}$, also called a norm matrix, such that

$$\mathbf{D}^- = \mathbf{M}^{-1}\left(\mathbf{Q}^- + \frac{1}{2}\mathbf{B}\right) \text{ and } \mathbf{D}^+ = \mathbf{M}^{-1}\left(\mathbf{Q}^+ + \frac{1}{2}\mathbf{B}\right), \tag{10}$$

where $\mathbf{B}$ is a boundary operator of the form $\mathbf{B} = \mathbf{t}_\beta \mathbf{t}_\beta^{\mathrm{T}} - \mathbf{t}_\alpha \mathbf{t}_\alpha^{\mathrm{T}}$, with boundary interpolation provided by $\mathbf{t}_\alpha^{\mathrm{T}} \mathbf{x}^l = \alpha^l$ and $\mathbf{t}_\beta^{\mathrm{T}} \mathbf{x}^l = \beta^l$ for $0 \leqslant l \leqslant r$ with $r \geqslant q$;

iii. integration by parts is mimicked by the property

$$\mathbf{Q}^+ + (\mathbf{Q}^-)^{\mathrm{T}} = \mathbf{0}; \tag{11}$$

iv. as an additional stability constraint, the symmetric matrix

$$\mathbf{C} := \frac{1}{2}(\mathbf{Q}^+ - \mathbf{Q}^-) = \frac{1}{2}\Big(\mathbf{Q}^+ + \big(\mathbf{Q}^+\big)^{\mathrm{T}}\Big) = -\frac{1}{2}\big(\mathbf{Q}^- + (\mathbf{Q}^-)^{\mathrm{T}}\big) \tag{12}$$

is negative semi-definite.

**Remark 2** Some remarks on the modifications with respect to classical SBP operators and on naming conventions are in order.

i. Due to the generalization provided by gSBP schemes beyond the traditional case of equidistant grid points including the boundary points $\alpha$ and $\beta$, the boundary operators of gSBP schemes extend the traditional boundary operator $\mathbf{B} = \mathrm{diag}(-1, 0, \cdots, 0, 1)$ to the given nodal set.
ii. In the case $\mathbf{C} = \mathbf{0}$, we have $\mathbf{D}^- = \mathbf{D}^+$ which is denoted a classical gSBP operator, since it does not feature any upwind directions.
iii. Upwind gSBP operators based on diagonal norm matrices $\mathbf{M}$ are denoted diagonal-norm upwind SBP operators.

The motivation for using upwind gSBP operators is the additional artificial damping introduced by enforcing the negative semi-definiteness of the matrix $\mathbf{C}$ in (12). We will illustrate this property of an upwind gSBP scheme in comparison with classical gSBP schemes in the following example.

Considering the spatial discretization of the linear advection equation (8) via an upwind gSBP scheme, we use the operator $\mathbf{D}^-$, whereas for $a < 0$, the operator $\mathbf{D}^+$ would be applied. We obtain

$$\frac{\mathrm{d}\mathbf{u}}{\mathrm{d}t} + a\mathbf{D}^-\mathbf{u} = \sigma\mathbf{M}^{-1}\mathbf{t}_\alpha\mathbf{t}_\alpha^{\mathrm{T}}\mathbf{u}, \tag{13}$$

where the right-hand side of the above equation contains the required simultaneous-approximation-term (SAT), dealing with the boundary condition at the left boundary. Thereby, the corresponding SAT parameter $\sigma \in \mathbb{R}$ is specified, such that a discrete energy estimate is fulfilled which mimics the continuous case, e.g., as in estimate (15).

Multiplying (13) from the left by $\mathbf{u}^{\mathrm{T}}\mathbf{M}$ and adding the transpose now yields

$$2\mathbf{u}^{\mathrm{T}}\mathbf{M}\frac{\mathrm{d}\mathbf{u}}{\mathrm{d}t} + a\mathbf{u}^{\mathrm{T}}\big(\mathbf{M}\,\mathbf{D}^- + (\mathbf{D}^-)^{\mathrm{T}}\mathbf{M}\big)\mathbf{u} = 2\sigma u_\alpha^2, \qquad u_\alpha = \mathbf{t}_\alpha^{\mathrm{T}}\mathbf{u}. \tag{14}$$

Using the upwind gSBP properties (11) and (12), we have

$$\mathbf{M}\,\mathbf{D}^- + (\mathbf{D}^-)^{\mathrm{T}}\mathbf{M} = \mathbf{Q}^- + (\mathbf{Q}^-)^{\mathrm{T}} + \mathbf{B} = -2\mathbf{C} + \mathbf{B},$$

and, therefore, setting $u_\beta = \mathbf{t}_\beta^{\mathrm{T}}\mathbf{u}$,

$$\frac{\mathrm{d}}{\mathrm{d}t}\|\mathbf{u}\|_{\mathbf{M}}^2 + a(u_\beta^2 - u_\alpha^2) = 2\sigma u_\alpha^2 + 2\mathbf{u}^{\mathrm{T}}\mathbf{C}\,\mathbf{u}.$$

For $\sigma \leqslant -\frac{a}{2}$, the time evolution of the discrete energy may now be estimated by

$$\frac{\mathrm{d}}{\mathrm{d}t}\|\mathbf{u}\|_{\mathbf{M}}^2 + a u_\beta^2 \leqslant 2\mathbf{u}^{\mathrm{T}}\mathbf{C}\,\mathbf{u} \leqslant 0, \tag{15}$$

where the term containing the negative semi-definite matrix $\mathbf{C}$ introduces additional artificial dissipation compared to the continuous estimate (9), and this additional dissipation vanishes for classical gSBP operators without upwind character, since $\mathbf{C} = \mathbf{0}$ in this case.

The upwind gSBP framework also enables the construction of second-derivative operators of the form

$$\begin{cases} \mathbf{D}_2^- = \mathbf{D}^+\mathbf{D}^- = \mathbf{M}^{-1}(-(\mathbf{D}^-)^{\mathrm{T}}\mathbf{M}\mathbf{D}^- + \mathbf{B}\mathbf{D}^-), \\ \mathbf{D}_2^+ = \mathbf{D}^-\mathbf{D}^+ = \mathbf{M}^{-1}(-(\mathbf{D}^+)^{\mathrm{T}}\mathbf{M}\mathbf{D}^+ + \mathbf{B}\mathbf{D}^+), \end{cases} \tag{16}$$

which approximate the second derivative $\frac{\partial^2}{\partial x^2}$.

### 2.2.1 Global DG Discretization of Advective Terms in Upwind SBP Framework

While the elementwise DG discretization of the linear advection equation results in a classical gSBP scheme on each grid cell, the corresponding global DG formulation over all grid cells results in its upwind variant. The classical SBP properties of element-wise 1D-DG schemes on general nodal sets are by now well-known, see, e.g., [15, 26, 28, 31]. In fact, the cellwise formulation (7) represents a gSBP scheme with the first-derivative gSBP operator $\widehat{\mathbf{D}} = \widehat{\mathbf{M}}^{-1}\mathbf{S}$ which approximates the continuous derivative $\frac{\partial}{\partial \xi}$ to degree $q = N$ and the degree of $\mathbf{B}$ is $r = q = N$, see [26, Lemma 1].

The boundary operator $\mathbf{B}$ is given by $\mathbf{B} = [\mathbf{L}\,\mathbf{L}^{\mathrm{T}}]_{-1}^{1}$, hence, for a single cell $(x_i, x_{i+1})$, we set $\alpha = x_i$, $\beta = x_{i+1}$ and obtain $\mathbf{t}_\alpha = \mathbf{L}(-1)$ and $\mathbf{t}_\beta = \mathbf{L}(1)$.

Next, we consider the global DG formulation over all grid cells. According to (7), on interior cells $I_i$, $i = 2, \cdots, K-1$, the cell-wise DG scheme in matrix-vector formulation applied to the linear advection equation (8), with the positive advection velocity $a > 0$, is given by

$$\frac{\Delta x_i}{2}\frac{\mathrm{d}\mathbf{u}^i}{\mathrm{d}t} + a\widehat{\mathbf{D}}\,\mathbf{u}^i = \widehat{\mathbf{M}}^{-1}\left[\left(au_h^i - (au)^{*,i}\right)\mathbf{L}\right]_{-1}^{1}, \quad i = 2, \cdots, K-1. \tag{17}$$

Hereby, a suitable numerical flux function $(au)^{*,i}$ is employed, where $(au)^{*,i}(-1) = (au)_i^*$ and $(au)^{*,i}(1) = (au)_{i+1}^*$ with

$$(au)_i^* = a\left(\frac{1}{2} + \theta\right)u_h^{i-1}(1) + a\left(\frac{1}{2} - \theta\right)u_h^i(-1), \quad \theta \in \left[0, \frac{1}{2}\right]. \tag{18}$$

Therefore, the choice of numerical flux functions ranges between the central flux for $\theta = 0$ and the upwind flux for $\theta = \frac{1}{2}$. Using the SAT boundary treatment with the SAT parameter $\sigma$, the scheme on the left boundary cell is given by

$$\frac{\Delta x_1}{2}\frac{\mathrm{d}\mathbf{u}^1}{\mathrm{d}t} + a\widehat{\mathbf{D}}\,\mathbf{u}^1 = \widehat{\mathbf{M}}^{-1}\left(au_h^1(1) - (au)_2^*\right)\mathbf{L}(1) + \sigma\widehat{\mathbf{M}}^{-1}\mathbf{L}(-1)\mathbf{L}(-1)^{\mathrm{T}}\mathbf{u}^1, \tag{19}$$

while on the rightmost DG cell with outgoing information through the right boundary, we have

$$\frac{\Delta x_K}{2}\frac{\mathrm{d}\mathbf{u}^K}{\mathrm{d}t} + a\widehat{\mathbf{D}}\,\mathbf{u}^K = -\widehat{\mathbf{M}}^{-1}\left(au_h^K(-1) - (au)_K^*\right)\mathbf{L}(-1). \tag{20}$$

We now rewrite the global DG discretization on the computational domain $\Omega$ as an upwind gSBP scheme which includes the interactions of degrees of freedom on adjacent cells. Under this assumption, we may rewrite (17) as

$$\frac{\mathrm{d}\mathbf{u}}{\mathrm{d}t} + a\mathbf{D}^{-}\mathbf{u} = \underline{\mathrm{SAT}}, \tag{21}$$

with the global nodal DG solution vector given by $\mathbf{u} = (\mathbf{u}^1, \cdots, \mathbf{u}^K)^{\mathrm{T}}$ and the extended SATs $\underline{\mathrm{SAT}} = \left(\underline{\mathrm{SAT}}_1^{\mathrm{T}}, 0, \cdots, 0\right)^{\mathrm{T}}$, where $\underline{\mathrm{SAT}}_1 = \frac{2\sigma}{\Delta x_1}\widehat{\mathbf{M}}^{-1}\mathbf{L}(-1)\mathbf{L}(-1)^{\mathrm{T}}\mathbf{u}^1$.

Now, the global DG derivative operator $\mathbf{D}^{-}$ is a block tridiagonal matrix with blocks corresponding to each DG element and its left and right adjacent cells. More precisely, we have

**Proposition 1** *The global DG derivative operator* $\mathbf{D}^{-} = \mathbf{D}^{-}(\theta)$ *in* (21) *is of the block tridiagonal form*

$$\mathbf{D}^{-}(\theta) = \begin{pmatrix} \frac{2}{\Delta x_1} & & & \\ & \frac{2}{\Delta x_2} & & \\ & & \ddots & \\ & & & \frac{2}{\Delta x_K} \end{pmatrix} \begin{pmatrix} \mathbf{A}_{\mathrm{lb}}(\theta) & \mathbf{A}_{12}(\theta) & & & \\ \mathbf{A}_{21}(\theta) & \mathbf{A}_{11}(\theta) & \mathbf{A}_{12}(\theta) & & \\ & \ddots & \ddots & \ddots & \\ & & \mathbf{A}_{21}(\theta) & \mathbf{A}_{11}(\theta) & \mathbf{A}_{12}(\theta) \\ & & & \mathbf{A}_{21}(\theta) & \mathbf{A}_{\mathrm{rb}}(\theta) \end{pmatrix},$$

*where the repeating blocks corresponding to interior cells and the left and right boundary blocks are given by*

$$\begin{aligned} \mathbf{A}_{11}(\theta) &= \widehat{\mathbf{D}} - \left(\frac{1}{2} - \theta\right)\widehat{\mathbf{M}}^{-1}\mathbf{L}(1)\mathbf{L}(1)^{\mathrm{T}} + \left(\frac{1}{2} + \theta\right)\widehat{\mathbf{M}}^{-1}\mathbf{L}(-1)\mathbf{L}(-1)^{\mathrm{T}}, \\ \mathbf{A}_{12}(\theta) &= \left(\frac{1}{2} - \theta\right)\widehat{\mathbf{M}}^{-1}\mathbf{L}(1)\mathbf{L}(-1)^{\mathrm{T}}, \\ \mathbf{A}_{21}(\theta) &= -\left(\frac{1}{2} + \theta\right)\widehat{\mathbf{M}}^{-1}\mathbf{L}(-1)\mathbf{L}(1)^{\mathrm{T}}, \\ \mathbf{A}_{\mathrm{lb}}(\theta) &= \widehat{\mathbf{D}} - \left(\frac{1}{2} - \theta\right)\widehat{\mathbf{M}}^{-1}\mathbf{L}(1)\mathbf{L}(1)^{\mathrm{T}}, \\ \mathbf{A}_{\mathrm{rb}}(\theta) &= \widehat{\mathbf{D}} + \left(\frac{1}{2} + \theta\right)\widehat{\mathbf{M}}^{-1}\mathbf{L}(-1)\mathbf{L}(-1)^{\mathrm{T}}, \end{aligned}$$

*where* $\widehat{\mathbf{D}}$ *denotes the cellwise first-derivative gSBP operator of the DG scheme* (7).

**Proof** The evaluation of the numerical flux function at an interface given in (18) may be rewritten as

$$(au)_i^* = a\left(\frac{1}{2} + \theta\right)\mathbf{L}(1)^{\mathrm{T}}\mathbf{u}^{i-1} + a\left(\frac{1}{2} - \theta\right)\mathbf{L}(-1)^{\mathrm{T}}\mathbf{u}^i$$

and the evaluation of the flux function $au_h^i$ at cell boundaries is obtained as

$$au_h^i(\pm 1) = a\mathbf{L}(\pm 1)^{\mathrm{T}}\mathbf{u}^i.$$

The desired block structure is obtained by inserting these contributions into (17) for interior DG cells and into (19) and (20) for boundary DG cells. Hereby, we account for the dependencies of the degrees of freedom on cell $i$ only on $\mathbf{u}^i$, $\mathbf{u}^{i-1}$, and $\mathbf{u}^{i+1}$ and set up the blocks accordingly.

The specific form of $\mathbf{D}^-$ allows us to show that the dual pair $\{\mathbf{D}^-, \mathbf{D}^+\}$ with $\mathbf{D}^- = \mathbf{D}^-(\theta)$ and a suitably constructed dual operator $\mathbf{D}^+$ satisfy the upwind gSBP properties specified in Definition 1 with $\alpha = x_a$ and $\beta = x_b$.

**Theorem 1** *The dual pair of discrete derivative operators*

$$\mathbf{D}^- = \mathbf{D}^-(\theta), \qquad \mathbf{D}^+ = \mathbf{D}^-(-\theta) \tag{22}$$

*is a dual pair of diagonal-norm upwind gSBP operators with respect to the global diagonal norm matrix*

$$\mathbf{M} = \operatorname{diag}\left(\frac{\Delta x_1}{2}\widehat{\mathbf{M}}, \cdots, \frac{\Delta x_K}{2}\widehat{\mathbf{M}}\right),$$

*and the generalized boundary operator*

$$\mathbf{B}_{\text{glob}} = \operatorname{diag}\left(\mathbf{B}_{\text{l}}, 0, \cdots, 0, \mathbf{B}_{\text{r}}\right)$$

*with* $\mathbf{B}_{\text{l}} = -\mathbf{L}(-1)\mathbf{L}(-1)^{\text{T}}$ *and* $\mathbf{B}_{\text{r}} = \mathbf{L}(1)\mathbf{L}(1)^{\text{T}}$.

**Proof** We define the matrices $\mathbf{Q}^-$ and $\mathbf{Q}^+$ in (10) by

$$\begin{aligned}
\mathbf{Q}^- &= \mathbf{M}\,\mathbf{D}^-(\theta) - \frac{1}{2}\mathbf{B}_{\text{glob}},\\
\mathbf{Q}^+ &= \mathbf{M}\,\mathbf{D}^+(\theta) - \frac{1}{2}\mathbf{B}_{\text{glob}}, \quad \mathbf{D}^+(\theta) = \mathbf{D}^-(-\theta).
\end{aligned}$$

Next, we show that $\mathbf{Q}^-$ and $\mathbf{Q}^+$ satisfy the SBP property (11), i.e., that the dual operator $\mathbf{D}^+$ is chosen properly. We have

$$\begin{aligned}
\mathbf{Q}^+ + (\mathbf{Q}^-)^{\text{T}} &= \mathbf{M}\,\mathbf{D}^-(-\theta) + (\mathbf{D}^-(\theta))^{\text{T}}\mathbf{M} - \mathbf{B}_{\text{glob}}\\
&= \begin{pmatrix}
\mathbf{Q}_{\text{lb}} & \mathbf{Q}_{12} & & & \\
\mathbf{Q}_{21} & \mathbf{Q}_{11} & \mathbf{Q}_{12} & & \\
 & \ddots & \ddots & \ddots & \\
 & & \mathbf{Q}_{21} & \mathbf{Q}_{11} & \mathbf{Q}_{12}\\
 & & & \mathbf{Q}_{21} & \mathbf{Q}_{\text{rb}}
\end{pmatrix},
\end{aligned}$$

where $\mathbf{Q}^+ + (\mathbf{Q}^-)^{\text{T}} = \mathbf{0}$, since the block representation given in Proposition 1 yields

$$\begin{aligned}
\mathbf{Q}_{11} &= \widehat{\mathbf{M}}\mathbf{A}_{11}(-\theta) + \mathbf{A}_{11}^{\text{T}}(\theta)\widehat{\mathbf{M}} = \widehat{\mathbf{M}}\widehat{\mathbf{D}} + \widehat{\mathbf{D}}^{\text{T}}\widehat{\mathbf{M}} - \mathbf{L}(1)\mathbf{L}(1)^{\text{T}} + \mathbf{L}(-1)\mathbf{L}(-1)^{\text{T}}\\
&= \widehat{\mathbf{M}}\widehat{\mathbf{D}} + \widehat{\mathbf{D}}^{\text{T}}\widehat{\mathbf{M}} - \mathbf{B} = \mathbf{0},\\
\mathbf{Q}_{\text{lb}} &= \widehat{\mathbf{M}}\mathbf{A}_{\text{lb}}(-\theta) + \mathbf{A}_{\text{lb}}^{\text{T}}(\theta)\widehat{\mathbf{M}} - \mathbf{B}_{\text{l}} = \widehat{\mathbf{M}}\widehat{\mathbf{D}} + \widehat{\mathbf{D}}^{\text{T}}\widehat{\mathbf{M}} - \mathbf{L}(1)\mathbf{L}(1)^{\text{T}} + \mathbf{L}(-1)\mathbf{L}(-1)^{\text{T}} = \mathbf{0},\\
\mathbf{Q}_{\text{rb}} &= \widehat{\mathbf{M}}\mathbf{A}_{\text{rb}}(-\theta) + \mathbf{A}_{\text{rb}}^{\text{T}}(\theta)\widehat{\mathbf{M}} - \mathbf{B}_{\text{r}} = \widehat{\mathbf{M}}\widehat{\mathbf{D}} + \widehat{\mathbf{D}}^{\text{T}}\widehat{\mathbf{M}} + \mathbf{L}(-1)\mathbf{L}(-1)^{\text{T}} - \mathbf{L}(1)\mathbf{L}(1)^{\text{T}} = \mathbf{0},\\
\mathbf{Q}_{12} &= \widehat{\mathbf{M}}\mathbf{A}_{12}(-\theta) + \mathbf{A}_{21}^{\text{T}}(\theta)\widehat{\mathbf{M}} = \left(\frac{1}{2} + \theta\right)\mathbf{L}(1)\mathbf{L}(-1)^{\text{T}} - \left(\frac{1}{2} + \theta\right)\mathbf{L}(1)\mathbf{L}(-1)^{\text{T}} = \mathbf{0},\\
\mathbf{Q}_{21} &= \widehat{\mathbf{M}}\mathbf{A}_{21}(-\theta) + \mathbf{A}_{12}^{\text{T}}(\theta)\widehat{\mathbf{M}} = -\left(\frac{1}{2} - \theta\right)\mathbf{L}(-1)\mathbf{L}(1)^{\text{T}} + \left(\frac{1}{2} - \theta\right)\mathbf{L}(-1)\mathbf{L}(1)^{\text{T}} = \mathbf{0}.
\end{aligned}$$

It remains to show the stability constraint (12), i.e., to show that

$$\mathbf{C} = \frac{1}{2}\left(\mathbf{Q}^{+} - \mathbf{Q}^{-}\right) = \frac{1}{2}\mathbf{M}(\mathbf{D}^{-}(-\theta) - \mathbf{D}^{-}(\theta)) \tag{23}$$

is negative semi-definite. We have

$$\mathbf{C} = \begin{pmatrix} \mathbf{C}_{\mathrm{lb}} & \mathbf{C}_{12} & & & \\ \mathbf{C}_{21} & \mathbf{C}_{11} & \mathbf{C}_{12} & & \\ & \ddots & \ddots & \ddots & \\ & & \mathbf{C}_{21} & \mathbf{C}_{11} & \mathbf{C}_{12} \\ & & & \mathbf{C}_{21} & \mathbf{C}_{\mathrm{rb}} \end{pmatrix}$$

with

$$\begin{aligned}
\mathbf{C}_{11} &= \frac{1}{2}\left(\widehat{\mathbf{M}}\mathbf{A}_{11}(-\theta) - \widehat{\mathbf{M}}\mathbf{A}_{11}(\theta)\right) = -\theta\left(\mathbf{L}(1)\mathbf{L}(1)^{\mathrm{T}} + \mathbf{L}(-1)\mathbf{L}(-1)^{\mathrm{T}}\right),\\
\mathbf{C}_{\mathrm{lb}} &= \frac{1}{2}\left(\widehat{\mathbf{M}}\mathbf{A}_{\mathrm{lb}}(-\theta) - \widehat{\mathbf{M}}\mathbf{A}_{\mathrm{lb}}(\theta)\right) = -\theta\mathbf{L}(1)\mathbf{L}(1)^{\mathrm{T}},\\
\mathbf{C}_{\mathrm{rb}} &= \frac{1}{2}\left(\widehat{\mathbf{M}}\mathbf{A}_{\mathrm{rb}}(-\theta) - \widehat{\mathbf{M}}\mathbf{A}_{\mathrm{rb}}(\theta)\right) = -\theta\mathbf{L}(-1)\mathbf{L}(-1)^{\mathrm{T}},\\
\mathbf{C}_{12} &= \frac{1}{2}\left(\widehat{\mathbf{M}}\mathbf{A}_{12}(-\theta) - \widehat{\mathbf{M}}\mathbf{A}_{12}(\theta)\right) = \theta\mathbf{L}(1)\mathbf{L}(-1)^{\mathrm{T}},\\
\mathbf{C}_{21} &= \frac{1}{2}\left(\widehat{\mathbf{M}}\mathbf{A}_{21}(-\theta) - \widehat{\mathbf{M}}\mathbf{A}_{21}(\theta)\right) = \theta\mathbf{L}(-1)\mathbf{L}(1)^{\mathrm{T}} = \mathbf{C}_{12}^{\mathrm{T}}.
\end{aligned}$$

Hence, we obtain

$$\begin{aligned}
\mathbf{u}^{\mathrm{T}}\mathbf{C}\,\mathbf{u} &= (\mathbf{u}^{1})^{\mathrm{T}}\mathbf{C}_{\mathrm{lb}}\,\mathbf{u}^{1} + \sum_{i=2}^{K-1}(\mathbf{u}^{i})^{\mathrm{T}}\mathbf{C}_{11}\,\mathbf{u}^{i} + (\mathbf{u}^{K})^{\mathrm{T}}\mathbf{C}_{\mathrm{rb}}\,\mathbf{u}^{K} + 2\sum_{i=1}^{K-1}(\mathbf{u}^{i})^{\mathrm{T}}\mathbf{C}_{12}\,\mathbf{u}^{i+1}\\
&= -\theta\left( (u_h^1(1))^2 + \sum_{i=2}^{K-1}\left((u_h^i(-1))^2 + (u_h^i(1))^2\right) + (u_h^K(-1))^2 - 2\sum_{i=1}^{K-1} u_h^i(1)u_h^{i+1}(-1)\right)\\
&= -\theta\sum_{i=1}^{K-1}\left(u_h^{i+1}(-1) - u_h^i(1)\right)^2 \leqslant 0,
\end{aligned}$$

i.e., $\mathbf{C}$ is negative semi-definite.

The following remark provides additional insight into the upwind gSBP character of the global DG formulation with numerical fluxes ranging from upwind to central fluxes.

**Remark 3** Since $\mathbf{D}^{+}(\theta) = \mathbf{D}^{-}(-\theta)$, we note that

i. the dual operator $\mathbf{D}^{+}$ arises naturally from the DG scheme applied to the linear advection equation with the negative advection velocity $a < 0$, and
ii. for $\theta = 0$, the global DG operator $\mathbf{D} = \mathbf{D}^{-}(0) = \mathbf{D}^{+}(0)$ is a classical diagonal-norm gSBP operator with global central character.

***Periodic first-derivative upwind gSBP operators***

Since the energy stability analysis in Sect. 3 regarding fully discrete IMEX upwind gSBP schemes is carried out for periodic advection-diffusion problems, the corresponding semi-discretization of (1) by upwind gSBP operators is devoid of SATs and periodic upwind gSBP operators are required which are of a slightly modified form, given by

$$\mathbf{D}^{-}(\theta)=\begin{pmatrix}\frac{2}{\Delta x_1} & & & \\ & \frac{2}{\Delta x_2} & & \\ & & \ddots & \\ & & & \frac{2}{\Delta x_K}\end{pmatrix}\begin{pmatrix}\mathbf{A}_{11}(\theta) & \mathbf{A}_{12}(\theta) & & \mathbf{A}_{21}(\theta)\\ \mathbf{A}_{21}(\theta) & \ddots & \ddots & \\ & \ddots & \ddots & \mathbf{A}_{12}(\theta)\\ \mathbf{A}_{12}(\theta) & & \mathbf{A}_{21}(\theta) & \mathbf{A}_{11}(\theta)\end{pmatrix}, \tag{24}$$

$$\mathbf{D}^{+}(\theta)=\mathbf{D}^{-}(-\theta), \tag{25}$$

$$\mathbf{D}_2(\theta)=\mathbf{D}^{-}(\theta)\mathbf{D}^{+}(\theta). \tag{26}$$

The construction of these periodic global discrete derivative operators is only based on the interior cell discretization (17). Therefore, they use the same building blocks given in Proposition 1 except for the absence of the boundary blocks $\mathbf{A}_{\mathrm{lb}}(\theta), \mathbf{A}_{\mathrm{rb}}(\theta)$.

### 2.2.2 Rewriting LDG and BR1 Schemes as Second-Derivative Upwind SBP Operators

The upwind gSBP operators resulting from the DG discretization of the linear advection equation in one space dimension can be used to construct specific second-derivative operators. In particular, the construction (16) of second-derivative upwind gSBP operators by combining $\mathbf{D}^{-}(\theta)$ and $\mathbf{D}^{+}(\theta)$ is directly related the DG technique of rewriting diffusion equations into a first-order system using an auxiliary variable for the solution derivative as in [3, 6, 7, 11]. For the linear heat equation

$$u_t = cu_{xx}, \quad (x,t)\in(x_\alpha, x_\beta)\times(0,T), \tag{27}$$

this procedure provides the system

$$u_t = cq_x, \quad q = u_x, \tag{28}$$

using an auxiliary variable $q$. For periodic boundary conditions, only the interior cell discretization (7) is relevant and the element-wise strong form DG semi-discretization yields

$$\begin{aligned}\frac{\Delta x_i}{2}\frac{\mathrm{d}\mathbf{u}^i}{\mathrm{d}t} &= c\widehat{\mathbf{D}}\mathbf{q}^i - c\widehat{\mathbf{M}}^{-1}\left[\left(q_h^i - q^{*,i}\right)\mathbf{L}\right]_{-1}^{1},\\ \frac{\Delta x_i}{2}\mathbf{q}^i &= \widehat{\mathbf{D}}\mathbf{u}^i - \widehat{\mathbf{M}}^{-1}\left[\left(u_h^i - u^{*,i}\right)\mathbf{L}\right]_{-1}^{1}.\end{aligned}$$

Regarding the specification of the numerical fluxes, at each element boundary, the left-hand side and right-hand side values of a piecewise continuous function $v$ are denoted by $v^-$ and $v^+$, respectively. The corresponding jump at element interfaces is denoted by $[v]=v^+-v^-$ and the arithmetic mean is given by $\{v\}=\frac{1}{2}(v^+ + v^-)$. The simplest choice for $q^*$ and $u^*$ is the BR1 scheme [6] given by the arithmetic means

$$q^{*,\mathrm{BR1}} = \{q\}, \quad u^{*,\mathrm{BR1}} = \{u\}. \tag{29}$$

Furthermore, the original LDG [11] yields a parameter-dependent family of diffusion fluxes

$$q^{*,\mathrm{LDG}} = \{q\} - c_{12}[q] + c_{11}[u], \quad u^{*,\mathrm{LDG}} = \{u\} + c_{12}[u], \tag{30}$$

which, in one space dimension, contains the BR1 approach as a specific case with $c_{11} = c_{12} = 0$. Commonly, LDG fluxes refer to the choice of alternating fluxes with $c_{12} = \pm 1$, $c_{11} = 0$, which offers two different variants. One implementation is thus given by

$$q^{*,\mathrm{LDG}_a} = q^-, \quad u^{*,\mathrm{LDG}_a} = u^+, \tag{31}$$

the second variant is specified by

$$q^{*,\mathrm{LDG}_b} = q^+, \quad u^{*,\mathrm{LDG}_b} = u^-. \tag{32}$$

For the $\mathrm{LDG}_a$, the $\mathrm{LDG}_b$, and the BR1 schemes, the numerical fluxes $q^*$ and $u^*$ can now be written as

$$\begin{aligned} q^* &= \left(\frac{1}{2} + \theta_q\right) q^- + \left(\frac{1}{2} - \theta_q\right) q^+, \\ u^* &= \left(\frac{1}{2} + \theta_u\right) u^- + \left(\frac{1}{2} - \theta_u\right) u^+, \end{aligned}$$

with $\theta_q = \frac{1}{2}, \theta_u = -\frac{1}{2}$ for $\mathrm{LDG}_a$, $\theta_q = -\frac{1}{2}, \theta_u = \frac{1}{2}$ for $\mathrm{LDG}_b$ and $\theta_q = \theta_u = 0$ for BR1.

Due to the assumed periodic boundary conditions, we now use the periodic global upwind gSBP operators (24), (25), and (26) developed in Sect. 2.2.1.

Thus, the $\mathrm{LDG}_a$ scheme can be rewritten as

$$\frac{\mathrm{d}\mathbf{u}}{\mathrm{d}t} = c\,\mathbf{D}^-\,\mathbf{q}, \quad \mathbf{q} = \mathbf{D}^+\,\mathbf{u},$$

where $\mathbf{u} = \left(\mathbf{u}^1, \cdots, \mathbf{u}^K\right)^{\mathrm{T}}$, $\mathbf{q} = \left(\mathbf{q}^1, \cdots, \mathbf{q}^K\right)^{\mathrm{T}}$ denote the global nodal solution vectors and the periodic upwind gSBP operator $\mathbf{D}^- = \mathbf{D}^-\left(\theta = \frac{1}{2}\right)$ is taken from (24). Hence, for the $\mathrm{LDG}_a$ scheme, we have

$$\mathrm{LDG}_a: \frac{\mathrm{d}\mathbf{u}}{\mathrm{d}t} = c\,\mathbf{D}_2^+\,\mathbf{u} = c\,\mathbf{D}^-\mathbf{D}^+\,\mathbf{u}.$$

Analogously, for the $\mathrm{LDG}_b$ variant we obtain

$$\mathrm{LDG}_b: \frac{\mathrm{d}\mathbf{u}}{\mathrm{d}t} = c\,\mathbf{D}_2^-\,\mathbf{u} = c\,\mathbf{D}^+\mathbf{D}^-\,\mathbf{u}.$$

Both of these formulations correspond to true upwind gSBP operators approximating the second derivative as defined in (16). In contrast, for the BR1 scheme, we have

$$\mathrm{BR1}: \frac{\mathrm{d}\mathbf{u}}{\mathrm{d}t} = c\,\mathbf{D}_2\,\mathbf{u} = c\,\mathbf{D}\,\mathbf{D}\,\mathbf{u},$$

where $\mathbf{D} = \mathbf{D}^{-}(\theta = 0) = \mathbf{D}^{+}(\theta = 0)$ is a classical first-derivative gSBP operator, as already mentioned in Remark 3. Thus, the BR1 scheme represents a second-derivative gSBP operator which is obtained by applying the first-derivative gSBP operator twice.

#### 2.2.3 Preliminary Estimates for First- and Second-Derivative Upwind gSBP Operators

The fully discrete stability analysis for IMEX upwind gSBP schemes in Sect. 3.1 and Sect. 3.2 will be restricted to periodic problems and is based on the discrete energy norm induced by the diagonal positive definite norm matrix $\mathbf{M}$. Defining the inner product

$$(\mathbf{u}, \mathbf{v})_{\mathbf{M}} = \mathbf{u}^{\mathrm{T}}\mathbf{M}\mathbf{v},$$

a discrete energy norm is given by

$$\|\mathbf{u}\|_{\mathbf{M}} = \sqrt{(\mathbf{u}, \mathbf{u})_{\mathbf{M}}} = \sqrt{\mathbf{u}^{\mathrm{T}}\mathbf{M}\mathbf{u}}.$$

Furthermore, the Cauchy-Schwarz inequality

$$(\mathbf{u}, \mathbf{v})_{\mathbf{M}} \leqslant \|\mathbf{u}\|_{\mathbf{M}}\|\mathbf{v}\|_{\mathbf{M}}$$

will be useful. As already mentioned, we assume periodic boundary conditions and the semi-discretization of the linear advection-diffusion equation (1) results in

$$\frac{\mathrm{d}\mathbf{u}}{\mathrm{d}t} + a\,\mathbf{D}^{-}\mathbf{u} = c\,\mathbf{D}_2\mathbf{u},$$

based on the periodic upwind gSBP operators given in (24) and (26).

Using this notation, we have the following estimates and relations between the operators $\mathbf{D}^{\pm}$ and $\mathbf{D}_2$, which will be required for the subsequent discrete energy stability analysis.

**Proposition 2** *For the periodic upwind gSBP operators defined in* (24), (25), *and* (26), *and arbitrary nodal values* $\mathbf{u}, \mathbf{v}$ *we have*

i. *dissipativity of* $-\mathbf{D}^{-}$, *i.e.*,

$$(\mathbf{D}^{-}\mathbf{u}, \mathbf{u})_{\mathbf{M}} \geqslant 0; \tag{33}$$

ii. *dissipativity of* $\mathbf{D}_2$, *more precisely*

$$\left(\mathbf{D}_2\mathbf{u}, \mathbf{v}\right)_{\mathbf{M}} = -\left(\mathbf{D}^{+}\mathbf{u}, \mathbf{D}^{+}\mathbf{v}\right)_{\mathbf{M}}, \tag{34}$$

$$\left(\mathbf{D}_2\mathbf{u}, \mathbf{u}\right)_{\mathbf{M}} = -\|\mathbf{D}^{+}\mathbf{u}\|_{\mathbf{M}}^2 \leqslant 0; \tag{35}$$

iii. *a Cauchy-Schwarz-type inequality*

$$|(\mathbf{D}^{-}\mathbf{u}, \mathbf{v})_{\mathbf{M}}| \leqslant \|\mathbf{u}\|_{\mathbf{M}}\|\mathbf{D}^{+}\mathbf{v}\|_{\mathbf{M}}. \tag{36}$$

**Proof** First, making use of the periodicity, the upwind gSBP property (11) reduces to the equation:

$$\mathbf{M}\mathbf{D}^{+} + (\mathbf{D}^{-})^{\mathrm{T}}\mathbf{M} = \mathbf{Q}^{+} + (\mathbf{Q}^{-})^{\mathrm{T}} = \mathbf{0},$$

leading to $(\mathbf{D}^{-})^{\mathrm{T}}\mathbf{M} = -\mathbf{M}\mathbf{D}^{+}$ which we will frequently use in the following derivations.

Due to the symmetry of the inner product, we have

$$(\mathbf{D}^{-}\mathbf{u}, \mathbf{u})_{\mathbf{M}} = \frac{1}{2}\big((\mathbf{D}^{-}\mathbf{u}, \mathbf{u})_{\mathbf{M}} + (\mathbf{u}, \mathbf{D}^{-}\mathbf{u})_{\mathbf{M}}\big),$$

and the first assertion is thus proven by

$$(\mathbf{D}^{-}\mathbf{u}, \mathbf{u})_{\mathbf{M}} = \frac{1}{2}\big(\mathbf{u}^{\mathrm{T}}(\mathbf{Q}^{-})^{\mathrm{T}}\mathbf{u} + \mathbf{u}^{\mathrm{T}}\mathbf{Q}^{-}\mathbf{u}\big) = -\frac{1}{2}\mathbf{u}^{\mathrm{T}}\big(\mathbf{Q}^{+} - \mathbf{Q}^{-}\big)\mathbf{u} = -\mathbf{u}^{\mathrm{T}}\mathbf{C}\mathbf{u} \geqslant 0.$$

Regarding the second assertion, we obtain (34) by

$$\begin{aligned}\big(\mathbf{D}_2\mathbf{u}, \mathbf{v}\big)_{\mathbf{M}} = \big(\mathbf{D}^{-}\mathbf{D}^{+}\mathbf{u}, \mathbf{v}\big)_{\mathbf{M}} &= \mathbf{u}^{\mathrm{T}}(\mathbf{D}^{+})^{\mathrm{T}}(\mathbf{D}^{-})^{\mathrm{T}}\mathbf{M}\mathbf{v} = -(\mathbf{D}^{+}\mathbf{u})^{\mathrm{T}}\mathbf{M}\mathbf{D}^{+}\mathbf{v} \\ &= -\big(\mathbf{D}^{+}\mathbf{u}, \mathbf{D}^{+}\mathbf{v}\big)_{\mathbf{M}},\end{aligned}$$

and (35) directly follows. For (36), we may rewrite

$$|(\mathbf{D}^{-}\mathbf{u}, \mathbf{v})_{\mathbf{M}}| = |\mathbf{u}^{\mathrm{T}}(\mathbf{D}^{-})^{\mathrm{T}}\mathbf{M}\mathbf{v}| = |\mathbf{u}^{\mathrm{T}}\mathbf{M}\mathbf{D}^{+}\mathbf{v}| = \big(\mathbf{u}, \mathbf{D}^{+}\mathbf{v}\big)_{\mathbf{M}}$$

and the third assertion follows from the Cauchy-Schwarz inequality for the inner product $(\cdot, \cdot)_{\mathbf{M}}$.

We will, furthermore, use the Young's inequality in the form

$$|ab| \leqslant \frac{a^2}{4C} + Cb^2 \tag{37}$$

for any constant $C > 0$ and any $a, b \in \mathbb{R}$.

## 3 IMEX Schemes for Advection-Diffusion Problems

Implicit-explicit one-step time integration methods for advection-diffusion equations are generally based on partitioned Runge-Kutta schemes applied to a partitioned system of ordinary differential equations of the form

$$\frac{\mathrm{d}\mathbf{u}}{\mathrm{d}t} = F_1(t, \mathbf{u}) + F_2(t, \mathbf{u}). \tag{38}$$

Hereby, different RK schemes, i.e., an explicit and an implicit one, are applied to the components $F_1$ and $F_2$ of the split right-hand side.

As in [4, 9] we consider $(s+1)$-stage explicit Runge-Kutta methods coupled with implicit $s$-stage DIRK schemes, where the absissae $\mathbf{c}^1, \mathbf{c}^2$ of the explicit and implicit schemes, respectively, fulfill $\mathbf{c}^1 = \begin{pmatrix} 0 \\ \mathbf{c}^2 \end{pmatrix}$. The DIRK scheme is then recast into an $(s+1)$-stage scheme as well by padding the first row and first column with zeros. We thus obtain the Butcher table

$$\begin{array}{c|cccccc|cccccc}
0 & 0 & & & & & & 0 & & & & \\
c_2 & a^1_{21} & 0 & & & & & 0 & a^2_{22} & & & \\
c_3 & a^1_{31} & a^1_{32} & 0 & & & & 0 & a^2_{32} & a^2_{33} & & \\
\vdots & \vdots & \vdots & \ddots & \ddots & & & \vdots & \vdots & \ddots & \ddots & \\
c_{s+1} & a^1_{s+1,1} & a^1_{s+1,2} & \cdots & a^1_{s+1,s} & 0 & & 0 & a^2_{s+1,2} & \cdots & a^2_{s+1,s} & a^2_{s+1,s+1} \\
\hline
 & b^1_1 & b^1_2 & b^1_2 & \cdots & b^1_{s+1} & & 0 & b^2_2 & b^2_3 & \cdots & b^2_{s+1}
\end{array}$$

Applied to the system of ODEs (38), the IMEX-RK scheme then has the form

$$\begin{cases}
\mathbf{u}^{(1)} = \mathbf{u}^n, \\
\mathbf{u}^{(i)} = \mathbf{u}^n + \Delta t \sum_{j=1}^{s+1} a^1_{ij} F_1\left(t^{n,j}, \mathbf{u}^{(j)}\right) + \Delta t \sum_{j=1}^{s+1} a^2_{ij} F_1\left(t^{n,j}, \mathbf{u}^{(j)}\right), \; i = 2, \cdots, s+1, \\
\mathbf{u}^{n+1} = \mathbf{u}^n + \Delta t \sum_{j=1}^{s+1} b^1_j F_1\left(t^{n,j}, \mathbf{u}^{(j)}\right) + \Delta t \sum_{j=1}^{s+1} b^2_j F_1\left(t^{n,j}, \mathbf{u}^{(j)}\right),
\end{cases} \tag{39}$$

where $\mathbf{u}^n, \mathbf{u}^{n+1}$ denote the approximations at times $t^n$ and $t^{n+1} = t^n + \Delta t$ and $\mathbf{u}^{(i)}$ are the intermediate stage values corresponding to the intermediate times $t^{n,i}$ given by $t^{n,i} = t^n + c_i \Delta t$ with $c_i = \sum_{j=1}^{s} a^1_{ij} = \sum_{j=1}^{s} a^2_{ij}$.

In [9], Calvo et al. derived necessary conditions on the coefficients of IMEX schemes of the above form to obtain a time step restriction of $\Delta t = \mathcal{O}(c/a^2)$ for advection-diffusion equations discretized by the Fourier method. These conditions are similar to L-stability and include $a^2_{s+1,j} = b^2_j$, $j = 2, \cdots, s+1$ as a requirement for the implicit scheme.

In this work, for the $L^2$-stability analysis of IMEX-DG schemes, we consider the following first- and second-order IMEX schemes also used in [38].

*First-order*:

$$\begin{array}{c|cc|cc}
0 & 0 & 0 & 0 & 0 \\
1 & 1 & 0 & 0 & 1 \\
\hline
 & 1 & 0 & 0 & 1
\end{array} \tag{40}$$

*Second-order*:

$$\begin{array}{c|ccc|ccc}
0 & 0 & 0 & 0 & 0 & 0 & 0 \\
\gamma & \gamma & 0 & 0 & 0 & \gamma & 0 \\
1 & \delta & 1-\delta & 0 & 0 & 1-\gamma & \gamma \\
\hline
 & \delta & 1-\delta & 0 & 0 & 1-\gamma & \gamma
\end{array} \tag{41}$$

with $\gamma = 1 - \frac{\sqrt{2}}{2}$ and $\delta = 1 - \frac{1}{2\gamma}$, fulfilling the identity

$$\gamma - \delta = 1.$$

In addition, some of the numerical experiments carried out in Sect. 4 utilize the third-order IMEX scheme taken from [9] which is given by

$$\begin{array}{c|cccc|cccc} 0 & 0 & 0 & 0 & 0 & 0 & 0 & 0 & 0 \\ \gamma & \gamma & 0 & 0 & 0 & 0 & \gamma & 0 & 0 \\ \frac{1+\gamma}{2} & \frac{1+\gamma}{2}-a_1 & a_1 & 0 & 0 & 0 & \frac{1-\gamma}{2} & \gamma & 0 \\ 1 & 0 & 1-a_2 & a_2 & 0 & 0 & b_1 & b_2 & \gamma \\ \hline & 0 & b_1 & b_2 & \gamma & 0 & b_1 & b_2 & \gamma \end{array} \tag{42}$$

with

$$\begin{aligned} \gamma \approx\ & 0.435\,866\,521\,508\,459 \text{ the middle root of } 6x^3-18x^2+9x-1, \\ b_1 =& -1.5\gamma^2+4\gamma-0.25, \\ b_2 =& 1.5\gamma^2-5\gamma+1.25, \\ a_1 =& -0.35, \\ a_2 =& \frac{1/3-2\gamma^2-2b_2a_1\gamma}{\gamma(1-\gamma)}. \end{aligned}$$

## 3.1 Stability Analysis for the First-Order IMEX Scheme

In principle, the energy stability analysis for IMEX upwind gSBP schemes is similar to the $L^2$-stability analysis carried out for LDG schemes in [38] and for the $(\sigma,\mu)$-family in [27]. The most significant difference of the current approach is the direct use of the discrete energy generated by the upwind gSBP norm matrix $\mathbf{M}$ instead of the usual $L^2$-norm of the approximate solution. This leads to a very transparent energy estimate which is devoid of finite-element-type inverse constants and holds for upwind gSBP schemes of arbitrary order in space. In particular, the provided time step restriction does not depend on the polynomial degree in space. For the applicability of the following stability analysis, compatibility of the first-derivative and second-derivative upwind gSBP operators in the sense that $\mathbf{D}_2=\mathbf{D}^-\mathbf{D}^+$ is crucial.

Discretizing the linear advection-diffusion equation (1) with $a>0$ by first- and second-derivative upwind gSBP operators $\mathbf{D}^-$ and $\mathbf{D}_2$ in space and the first-order IMEX scheme (40) as time integrator yields the fully discrete IMEX upwind gSBP scheme

$$\mathbf{u}^{n+1}=\mathbf{u}^n-a\Delta t\,\mathbf{D}^-\mathbf{u}^n+c\Delta t\,\mathbf{D}_2\mathbf{u}^{n+1}. \tag{43}$$

Applying $(\cdot,\cdot)_{\mathbf{M}}$ on both sides of (43) with $\mathbf{u}^{n+1}$ as the second argument, we get

$$\left(\mathbf{u}^{n+1},\mathbf{u}^{n+1}\right)_{\mathbf{M}}=\left(\mathbf{u}^n,\mathbf{u}^{n+1}\right)_{\mathbf{M}}-a\Delta t\left(\mathbf{D}^-\mathbf{u}^n,\mathbf{u}^{n+1}\right)_{\mathbf{M}}+c\Delta t\left(\mathbf{D}_2\mathbf{u}^{n+1},\mathbf{u}^{n+1}\right)_{\mathbf{M}}. \tag{44}$$

Denoting $\mathbf{u}^\delta=\mathbf{u}^{n+1}-\mathbf{u}^n$ and using identities analogously to the derivation in [38] combined with the estimate (33), we obtain

$$\begin{aligned} \left(\mathbf{u}^{n+1}-\mathbf{u}^n,\mathbf{u}^{n+1}\right)_{\mathbf{M}} &= \frac{1}{2}\left(\mathbf{u}^{n+1},\mathbf{u}^{n+1}\right)_{\mathbf{M}}+\frac{1}{2}\left(\mathbf{u}^\delta,\mathbf{u}^\delta\right)_{\mathbf{M}}-\frac{1}{2}(\mathbf{u}^n,\mathbf{u}^n)_{\mathbf{M}}, \\ -\left(\mathbf{D}^-\mathbf{u}^n,\mathbf{u}^{n+1}\right)_{\mathbf{M}} &= \left(\mathbf{D}^-\mathbf{u}^\delta,\mathbf{u}^{n+1}\right)_{\mathbf{M}}-\left(\mathbf{D}^-\mathbf{u}^{n+1},\mathbf{u}^{n+1}\right)_{\mathbf{M}}\leqslant\left(\mathbf{D}^-\mathbf{u}^\delta,\mathbf{u}^{n+1}\right)_{\mathbf{M}}. \end{aligned}$$

Inserting this into (44) and employing (35), (36), and finally the Young's inequality (37), we have

$$\begin{aligned}\frac{1}{2}\left(\|\mathbf{u}^{n+1}\|_{\mathbf{M}}^2-\|\mathbf{u}^n\|_{\mathbf{M}}^2\right) &\leqslant -\frac{1}{2}\|\mathbf{u}^\delta\|_{\mathbf{M}}^2+a\Delta t\big(\mathbf{D}^-\mathbf{u}^\delta,\mathbf{u}^{n+1}\big)_{\mathbf{M}}+c\Delta t\big(\mathbf{D}_2\mathbf{u}^{n+1},\mathbf{u}^{n+1}\big)_{\mathbf{M}}\\ &\leqslant -\frac{1}{2}\|\mathbf{u}^\delta\|_{\mathbf{M}}^2+a\Delta t\|\mathbf{u}^\delta\|_{\mathbf{M}}\|\mathbf{D}^+\mathbf{u}^{n+1}\|_{\mathbf{M}}-c\Delta t\|\mathbf{D}^+\mathbf{u}^{n+1}\|_{\mathbf{M}}^2\\ &\leqslant -\frac{1}{2}\|\mathbf{u}^\delta\|_{\mathbf{M}}^2+\frac{a^2\Delta t}{4c}\|\mathbf{u}^\delta\|_{\mathbf{M}}^2+c\Delta t\|\mathbf{D}^+\mathbf{u}^{n+1}\|_{\mathbf{M}}^2-c\Delta t\|\mathbf{D}^+\mathbf{u}^{n+1}\|_{\mathbf{M}}^2\\ &=\left(\frac{a^2\Delta t}{4c}-\frac{1}{2}\right)\|\mathbf{u}^\delta\|_{\mathbf{M}}^2.\end{aligned}$$

Summarizing the above findings we have the following theorem.

**Theorem 2** *If* $\{\mathbf{D}^-,\mathbf{D}^+\}$ *is any dual pair of first-derivative upwind SBP operators with the norm matrix* $\mathbf{M}$ *and the second-derivative upwind SBP operator is chosen as* $\mathbf{D}_2=\mathbf{D}^-\mathbf{D}^+$, *then the fully discrete solution of the IMEX-upwind SBP scheme* (43) *satisfies the discrete energy stability of the form* $\|\mathbf{u}^{n+1}\|_{\mathbf{M}}\leqslant\|\mathbf{u}^n\|_{\mathbf{M}}$ *if the time step is bounded by*

$$\Delta t\leqslant\frac{2c}{a^2}.$$

It is worth noting that the assertion of Theorem 2 is quite general, since it also holds for the LDG diffusion scheme coupled with upwind fluxes for advection, i.e., the scheme analyzed in [38]. In addition, the constant in the above time step restriction is specifically determined and is independent of the polynomial degree of the space discretization, since no detour through finite-element spaces for $u$ and $u_x$ is required.

### 3.2 Stability Analysis for the Second-Order IMEX Scheme

Discretizing the linear advection-diffusion equation (1) by upwind gSBP operators in space as in Sect. 3.1 and using the second-order IMEX time integrator (41) yields the fully discrete IMEX upwind gSBP scheme

$$\mathbf{u}^{(n,1)}=\mathbf{u}^n+\Delta t\gamma\left(-a\,\mathbf{D}^-\mathbf{u}^n+c\,\mathbf{D}_2\mathbf{u}^{(n,1)}\right), \tag{45}$$

$$\begin{aligned}\mathbf{u}^{n+1}=&\mathbf{u}^n-\Delta t\big(\delta a\,\mathbf{D}^-\mathbf{u}^n+(1-\delta)a\,\mathbf{D}^-\mathbf{u}^{(n,1)}\big)\\ &+\Delta t\big((1-\gamma)c\,\mathbf{D}_2\mathbf{u}^{(n,1)}+\gamma c\,\mathbf{D}_2\mathbf{u}^{n+1}\big).\end{aligned} \tag{46}$$

Employing the discrete inner product $(\cdot,\cdot)_{\mathbf{M}}$ on both sides of (45) with $\mathbf{u}^{(n,1)}$ as the second argument and on both sides of (46) with $\mathbf{u}^{n+1}$ as the second argument, we have

$$(\mathbf{u}^{(n,1)}-\mathbf{u}^n,\mathbf{u}^{(n,1)})_{\mathbf{M}}=-a\gamma\Delta t\big(\mathbf{D}^-\mathbf{u}^n,\mathbf{u}^{(n,1)}\big)_{\mathbf{M}}+c\gamma\Delta t\big(\mathbf{D}_2\mathbf{u}^{(n,1)},\mathbf{u}^{(n,1)}\big)_{\mathbf{M}}, \tag{47}$$

$$\begin{aligned}(\mathbf{u}^{n+1}-\mathbf{u}^{(n,1)},\mathbf{u}^{n+1})_{\mathbf{M}}=&\,a\Delta t\Big(\big(\mathbf{D}^-\mathbf{u}^n,\mathbf{u}^{n+1}\big)_{\mathbf{M}}-(1-\delta)\big(\mathbf{D}^-\mathbf{u}^{(n,1)},\mathbf{u}^{n+1}\big)_{\mathbf{M}}\Big)\\ &+c\Delta t\Big((1-2\gamma)\big(\mathbf{D}_2\mathbf{u}^{(n,1)},\mathbf{u}^{n+1}\big)_{\mathbf{M}}+\gamma\big(\mathbf{D}_2\mathbf{u}^{n+1},\mathbf{u}^{n+1}\big)_{\mathbf{M}}\Big),\end{aligned} \tag{48}$$

where we have used $\gamma-\delta=1$. Setting $\mathbf{u}^{\delta,1}=\mathbf{u}^{(n,1)}-\mathbf{u}^n$ and $\mathbf{u}^{\delta,2}=\mathbf{u}^{n+1}-\mathbf{u}^{(n,1)}$ and adding (47) and (48), we obtain

$$\frac{1}{2}\left(\|\mathbf{u}^{n+1}\|^2 - \|\mathbf{u}^n\|^2 + \|\mathbf{u}^{\delta,1}\|^2 + \|\mathbf{u}^{\delta,2}\|^2\right) = \Delta t(R_1 + R_2) \tag{49}$$

with

$$\begin{aligned}R_1 &= -a\gamma\big(\mathbf{D}^-\mathbf{u}^n, \mathbf{u}^{(n,1)}\big)_{\mathbf{M}} + a\big(\mathbf{D}^-\mathbf{u}^n, \mathbf{u}^{n+1}\big)_{\mathbf{M}} - a(1-\delta)\big(\mathbf{D}^-\mathbf{u}^{(n,1)}, \mathbf{u}^{n+1}\big)_{\mathbf{M}}\\ &= -a\gamma\Big(\big(\mathbf{D}^-\mathbf{u}^{(n,1)}, \mathbf{u}^{(n,1)}\big)_{\mathbf{M}} - \big(\mathbf{D}^-\mathbf{u}^{\delta,1}, \mathbf{u}^{(n,1)}\big)_{\mathbf{M}}\Big) - a(1-\gamma)\big(\mathbf{D}^-\mathbf{u}^{n+1}, \mathbf{u}^{n+1}\big)_{\mathbf{M}}\\ &\quad - a\big(\mathbf{D}^-\mathbf{u}^{\delta,1}, \mathbf{u}^{n+1}\big)_{\mathbf{M}} + a(1-\gamma)\big(\mathbf{D}^-\mathbf{u}^{\delta,2}, \mathbf{u}^{n+1}\big)_{\mathbf{M}},\end{aligned}$$

and via (35),

$$R_2 = -c\big(\gamma\|\mathbf{D}^+\mathbf{u}^{(n,1)}\|_{\mathbf{M}}^2 + (1-2\gamma)\big(\mathbf{D}^+\mathbf{u}^{(n,1)}, \mathbf{D}^+\mathbf{u}^{n+1}\big)_{\mathbf{M}} + \gamma\|\mathbf{D}^+\mathbf{u}^{n+1}\|_{\mathbf{M}}^2\big).$$

Furthermore, $R_2$ fulfills

$$R_2 = S - \frac{1}{4}c\gamma\big(\|\mathbf{D}^+\mathbf{u}^{(n,1)}\|_{\mathbf{M}}^2 + \|\mathbf{D}^+\mathbf{u}^{n+1}\|_{\mathbf{M}}^2\big)$$

with

$$S = -c\Big(\frac{3}{4}\gamma\|\mathbf{D}^+\mathbf{u}^{(n,1)}\|_{\mathbf{M}}^2 + (1-2\gamma)\big(\mathbf{D}^+\mathbf{u}^{(n,1)}, \mathbf{D}^+\mathbf{u}^{n+1}\big)_{\mathbf{M}} + \frac{3}{4}\gamma\|\mathbf{D}^+\mathbf{u}^{n+1}\|_{\mathbf{M}}^2\Big) \leqslant 0$$

since $\mathbf{S} = \begin{pmatrix} \frac{3}{4}\gamma & \frac{1}{2}-\gamma \\ \frac{1}{2}-\gamma & \frac{3}{4}\gamma \end{pmatrix}$ can be proven positive definite for $\gamma = 1 - \frac{\sqrt{2}}{2}$ and we have

$$S = -c\,\mathbf{w}^{\mathrm{T}}(\mathbf{S}\otimes\mathbf{M})\mathbf{w} \quad \text{with } \mathbf{w} = \begin{pmatrix} \mathbf{D}^+\mathbf{u}^{(n,1)} \\ \mathbf{D}^+\mathbf{u}^{n+1} \end{pmatrix}.$$

Using (33), (36), we can, furthermore, bound $R_1$ by

$$\begin{aligned}R_1 &\leqslant a\big(\|\mathbf{u}^{\delta,1}\|_{\mathbf{M}}\big(\gamma\|\mathbf{D}^+\mathbf{u}^{(n,1)}\|_{\mathbf{M}} + \|\mathbf{D}^+\mathbf{u}^{n+1}\|_{\mathbf{M}}\big) + (1-\gamma)\|\mathbf{u}^{\delta,2}\|_{\mathbf{M}}\|\mathbf{D}^+\mathbf{u}^{n+1}\|_{\mathbf{M}}\big)\\ &\leqslant \frac{a^2}{c}\left(\left(\gamma + \frac{3}{2\gamma}\right)\|\mathbf{u}^{\delta,1}\|_{\mathbf{M}}^2 + \frac{11}{2}\|\mathbf{u}^{\delta,2}\|_{\mathbf{M}}^2\right) + c\frac{\gamma}{4}\|\mathbf{D}^+\mathbf{u}^{(n,1)}\|_{\mathbf{M}}^2\\ &\quad + c\left(\frac{\gamma}{6} + \frac{(1-\gamma)^2}{22}\right)\|\mathbf{D}^+\mathbf{u}^{n+1}\|_{\mathbf{M}}^2\\ &\leqslant \frac{11}{2}\frac{a^2}{c}\big(\|\mathbf{u}^{\delta,1}\|_{\mathbf{M}}^2 + \|\mathbf{u}^{\delta,2}\|_{\mathbf{M}}^2\big) + \frac{1}{4}\gamma c\big(\|\mathbf{D}^+\mathbf{u}^{(n,1)}\|_{\mathbf{M}}^2 + \|\mathbf{D}^+\mathbf{u}^{n+1}\|_{\mathbf{M}}^2\big),\end{aligned}$$

since for the parameter $\gamma = 1 - \frac{\sqrt{2}}{2}$, the inequalities $\gamma + \frac{3}{2\gamma} \leqslant \frac{11}{2}$ and $\frac{\gamma}{6} + \frac{(1-\gamma)^2}{22} \leqslant \frac{1}{4}$ hold.

Inserting the estimates for $R_1$ and $R_2$ into (49), we obtain the following theorem.

**Theorem 3** *If* $\{\mathbf{D}^-, \mathbf{D}^+\}$ *is any dual pair of first-derivative upwind SBP operators with the norm matrix M and the second-derivative upwind SBP operator is chosen as* $\mathbf{D}_2 = \mathbf{D}^-\mathbf{D}^+$, *then the fully discrete solution of the IMEX-upwind SBP schemes* (45), (46) *satisfies the discrete energy stability of the form* $\|\mathbf{u}^{n+1}\|_{\mathbf{M}} \leqslant \|\mathbf{u}^n\|_{\mathbf{M}}$ *if the time step is bounded by*

$$\Delta t \leqslant \frac{c}{11a^2}.$$

## 4 Numerical Results

In this section, the stability and accuracy of IMEX upwind gSBP schemes are compared with respect to different combinations of the first-derivative and second-derivative upwind gSBP operators. All upwind gSBP operators in this section are based on uniform grids and the Legendre-Gauss-Lobatto quadrature rule.

**Table 1** Stability of IMEX upwind gSBP schemes (43), (45), and (46) depending on $\left(\theta_{\mathtt{adv}}, \theta_{\mathtt{diff}}\right)$. Values of $\tau = \frac{a^2}{c}\Delta t_{\max}$, where $\Delta t_{\max}$ is the maximum time step ensuring a non-increasing discrete energy

| $N$ | $K$ | $a=c=0.1$ | | | | $a=0.2,\ c=0.01$ | | | |
|---|---|---|---|---|---|---|---|---|---|
| | | Upwind gSBP discretization $\left(\theta_{\mathtt{adv}}, \theta_{\mathtt{diff}}\right)$ | | | | | | | |
| | | $\left(\frac{1}{2},\frac{1}{2}\right)$ | $\left(\frac{1}{2},0\right)$ | $\left(\frac{1}{4},\frac{1}{4}\right)$ | $(0,0)$ | $\left(\frac{1}{2},\frac{1}{2}\right)$ | $\left(\frac{1}{2},0\right)$ | $\left(\frac{1}{4},\frac{1}{4}\right)$ | $(0,0)$ |
| First-order IMEX time integration | | | | | | | | | |
| 1 | 20 | + | 3.3E−01 | + | + | 2.1 | 2.0 | 2.1 | 2.0 |
| | 40 | + | 1.6E−01 | + | + | 2.0 | 2.0 | 2.0 | 2.0 |
| | 80 | + | 7.9E−02 | + | + | 2.0 | 1.6 | 2.0 | 2.0 |
| | 160 | + | 3.9E−02 | + | + | 2.0 | 7.9E−01 | 2.0 | 2.0 |
| | 320 | + | 2.0E−02 | + | + | 2.0 | 3.9E−01 | 2.0 | 2.0 |
| 2 | 20 | + | 1.5E−01 | + | + | 2.0 | 2.0 | 2.0 | 2.0 |
| | 40 | + | 7.8E−02 | + | + | 2.0 | 1.3 | 2.0 | 2.0 |
| | 80 | + | 3.9E−02 | + | + | 2.0 | 6.7E−01 | 2.0 | 2.0 |
| | 160 | + | 1.9E−02 | + | + | 2.0 | 3.5E−01 | 2.0 | 2.0 |
| | 320 | + | 9.8E−03 | + | + | 2.0 | 1.8e−01 | 2.0 | 2.0 |
| 3 | 20 | + | 1.0E−01 | + | + | 2.0 | 1.5 | 2.0 | 2.0 |
| | 40 | + | 5.2E−02 | + | + | 2.0 | 8.3E−01 | 2.0 | 2.0 |
| | 80 | + | 2.6E−02 | + | + | 2.0 | 4.6E−01 | 2.0 | 2.0 |
| | 160 | + | 1.3E−02 | + | + | 2.0 | 2.4E−01 | 2.0 | 2.0 |
| | 320 | + | 6.5E−03 | + | + | 2.0 | 1.2E−01 | 2.0 | 2.0 |
| Second-order IMEX time integration | | | | | | | | | |
| 1 | 20 | 2.4 | 3.2E−01 | 1.5 | 2.4 | 6.6 | 6.3 | 6.5 | 5.1 |
| | 40 | 2.4 | 1.6E−01 | 2.4 | 2.4 | 3.5 | 3.1 | 3.7 | 3.3 |
| | 80 | 2.4 | 7.9E−02 | 2.4 | 2.4 | 2.0 | 1.6 | 2.1 | 2.1 |
| | 160 | 2.4 | 3.9E−02 | 2.4 | 2.4 | 1.4 | 7.8E−01 | 1.5 | 1.4 |
| | 320 | 2.4 | 2.0E−02 | 2.4 | 2.4 | 1.5 | 3.9E−01 | 1.4 | 1.4 |
| 2 | 20 | 2.4 | 1.3E−01 | 2.4 | 2.4 | 2.9 | 2.3 | 3.1 | 2.8 |
| | 40 | 2.4 | 7.1E−02 | 2.4 | 2.4 | 1.9 | 1.2 | 1.9 | 1.8 |
| | 80 | 2.4 | 3.7E−02 | 2.4 | 2.4 | 1.7 | 6.0E−01 | 1.5 | 1.4 |
| | 160 | 2.4 | 1.9E−02 | 2.4 | 2.4 | 1.4 | 3.1E−01 | 1.5 | 1.4 |
| | 320 | 2.4 | 9.6E−03 | 2.4 | 2.4 | 1.4 | 1.6E−01 | 1.4 | 1.4 |
| 3 | 20 | 2.4 | 9.1E−02 | 2.4 | 2.4 | 2.2 | 1.3 | 2.1 | 1.9 |
| | 40 | 2.4 | 4.8E−02 | 2.4 | 2.4 | 1.8 | 6.8E−01 | 1.6 | 1.4 |
| | 80 | 2.4 | 2.5E−02 | 2.4 | 2.4 | 1.4 | 3.6E−01 | 1.4 | 1.4 |
| | 160 | 2.4 | 1.3E−02 | 2.4 | 2.4 | 1.4 | 2.0E−01 | 1.4 | 1.4 |
| | 320 | 2.4 | 6.5E−02 | 2.4 | 2.4 | 1.4 | 1.1E−01 | 1.4 | 1.4 |

### 4.1 Numerical Study on Time Step Restrictions

We consider the exact solution

$$u(x,t) = \mathrm{e}^{-ct}\sin(x - at) \tag{50}$$

of the linear advection-diffusion equation (1) in the interval $(x_a, x_b) = (-\pi, \pi)$. The advection-diffusion problem is discretized in space by different upwind gSBP operators on uniform grids using a number of $N + 1$ Legendre-Gauss-Lobatto nodes on each grid cell. The implemented upwind gSBP operators are specified by the parameters $\theta_{\mathtt{adv}}$ and $\theta_{\mathtt{diff}}$, respectively, where we choose $\mathbf{D}^- = \mathbf{D}^-(\theta = \theta_{\mathtt{adv}})$ and $\mathbf{D}_2 = \mathbf{D}^-(\theta_{\mathtt{diff}})\mathbf{D}^+(\theta_{\mathtt{diff}})$. Thus, for $\theta_{\mathtt{adv}} = \theta_{\mathtt{diff}}$, the first-derivative and second-derivative operators are compatible in the sense of Theorems 2 and 3. For instance, the choice $\left(\theta_{\mathtt{adv}}, \theta_{\mathtt{diff}}\right) = \left(\frac{1}{2}, \frac{1}{2}\right)$ represents upwind advection fluxes paired with LDG diffusion fluxes, while $\left(\theta_{\mathtt{adv}}, \theta_{\mathtt{diff}}\right) = (0, 0)$ denotes the pairing of central advection fluxes with the BR1 diffusion scheme. The first-order IMEX scheme (40) and the second-order IMEX scheme (41) are then used to discretize advection terms explicitly and diffusion terms implicitly in time. From the theoretical analysis, we expect the schemes to be stable for time steps $\Delta t \leqslant \tau \frac{c}{a^2}$, where $\tau$ is some constant independent of grid refinement.

For different advection and diffusion parameters $a$ and $c$, Table 1 lists the maximum stable time steps obtained by the IMEX schemes of first- and second-order coupled with upwind gSBP operators using $N \leqslant 3$. In the numerical setup, we successively double the number of cells $K$ and compute the numerical solution until the final time $T = 1\,000$ is reached. The maximum stable time step $\Delta t_{\max}$ is determined as the maximum time step for which the discrete energy norm $(\mathbf{u}, \mathbf{u})_{\mathbf{M}}$ of the numerical solution is non-increasing and the corresponding values of $\tau = \frac{a^2}{c}\Delta t_{\max}$ are listed. An entry of + indicates that arbitrarily large time steps may be chosen without increasing the discrete energy in time.

**Table 2** Values of $\tau = \frac{a^2}{c}\Delta t_{\max}$ for IMEX upwind gSBP schemes using the third-order IMEX method (42)

| $N$ | $K$ | $a = c = 0.1$ | | | | $a = 0.2,\ c = 0.01$ | | | |
|---|---|---|---|---|---|---|---|---|---|
| | | Upwind gSBP discretization $\left(\theta_{\mathtt{adv}}, \theta_{\mathtt{diff}}\right)$ | | | | | | | |
| | | $\left(\frac{1}{2}, \frac{1}{2}\right)$ | $\left(\frac{1}{2}, 0\right)$ | $\left(\frac{1}{4}, \frac{1}{4}\right)$ | $(0, 0)$ | $\left(\frac{1}{2}, \frac{1}{2}\right)$ | $\left(\frac{1}{2}, 0\right)$ | $\left(\frac{1}{4}, \frac{1}{4}\right)$ | $(0, 0)$ |
| 1 | 20 | 5.9 | 3.6E−01 | 5.9 | 5.9 | 7.5 | 7.0 | 10.1 | 11.1 |
| | 40 | 5.9 | 1.8E−01 | 5.9 | 5.9 | 4.1 | 3.5 | 6.4 | 7.2 |
| | 80 | 5.9 | 8.7E−02 | 5.9 | 5.9 | 2.6 | 1.7 | 3.9 | 5.5 |
| | 160 | 5.9 | 4.3E−02 | 5.9 | 5.9 | 2.4 | 8.6E−01 | 2.6 | 5.1 |
| | 320 | 5.9 | 2.2E−02 | 5.9 | 5.9 | 5.2 | 4.3E−01 | 2.6 | 5.0 |
| 2 | 20 | 5.9 | 1.4E−01 | 5.9 | 5.9 | 3.4 | 2.6 | 5.0 | 6.4 |
| | 40 | 5.9 | 6.9E−02 | 5.9 | 5.9 | 2.1 | 1.3 | 3.1 | 5.2 |
| | 80 | 5.9 | 3.5E−02 | 5.9 | 5.9 | 2.1 | 6.5E−01 | 2.4 | 5.0 |
| | 160 | 5.9 | 1.8E−02 | 5.9 | 5.9 | 4.6 | 3.3E−01 | 2.7 | 4.9 |
| | 320 | 5.9 | 8.9E−03 | 5.9 | 5.9 | 5.0 | 1.7E−01 | 5.0 | 4.9 |
| 3 | 20 | 5.9 | 8.1E−02 | 5.9 | 5.9 | 2.3 | 1.4 | 3.2 | 5.4 |
| | 40 | 5.9 | 4.1E−02 | 5.9 | 5.9 | 1.8 | 7.3E−01 | 2.5 | 5.0 |
| | 80 | 5.9 | 2.1E−02 | 5.9 | 5.9 | 3.9 | 3.7E−01 | 2.2 | 5.0 |
| | 160 | 5.9 | 1.0E−02 | 5.9 | 5.9 | 4.8 | 1.9E−01 | 4.2 | 5.0 |
| | 320 | 5.9 | 5.2E−03 | 5.9 | 5.9 | 5.0 | 9.9E−02 | 5.0 | 5.0 |

As predicted by the theoretical analysis, for all schemes using compatible first-derivative and second-derivative upwind gSBP operators with $\theta_{\mathtt{adv}} = \theta_{\mathtt{diff}}$, the maximum admissible time step is bounded from below under grid refinement. The precise bound varies for different advection-diffusion parameters $a$, $c$ but fulfills the theoretically predicted time step restriction as long as the compatibility condition is fulfilled. Moreover, the different schemes possess nearly the same time step bounds under the same setup with respect to the parameters $a$, $c$ and the chosen IMEX scheme. Regarding the enhanced IMEX stability property, the time step bound is lower for second-order IMEX time integration as predicted by the theoretical analysis.

On the other hand, as already noted in [27], the choice of $\left(\theta_{\mathtt{adv}}, \theta_{\mathtt{diff}}\right) = \left(\frac{1}{2}, 0\right)$ representing the BR1 diffusion scheme paired upwind advection fluxes does not admit grid-independent time step sizes, in accordance with the theoretical analysis with only covers the case $\theta_{\mathtt{adv}} = \theta_{\mathtt{diff}}$ as already stated. Instead, we observe admissible time steps of $\tau = \mathcal{O}(\Delta x)$, analogous to the time step restrictions for explicitly discretized advection equations. Replacing the upwind advection fluxes by central fluxes is compatible with the BR1 scheme and yields the enhanced IMEX stability which can also be observed in the results of Table 1.

Regarding higher order time discretization, Table 2 lists the corresponding maximum time steps $\tau$ ensuring a non-increasing discrete energy when combining the previously considered spatial upwind gSBP schemes for $N \leqslant 3$ with the third-order IMEX scheme (42). Although the theoretical analysis in Sect. 3 only covers the first-order and second-order IMEX schemes (40) and (41), respectively, an analogous stability behavior can be observed for the third-order IMEX scheme. In addition, for this setup, compatible first- and second-derivative upwind gSBP operators yield the enhanced IMEX stability, while for the combination of upwind advection fluxes with BR1 diffusion, i.e., $\left(\theta_{\mathtt{adv}}, \theta_{\mathtt{diff}}\right) = \left(\frac{1}{2}, 0\right)$, the classical CFL condition for the pure advection equation needs to be enforced.

### 4.2 Numerical Study on Accuracy of IMEX Upwind gSBP Schemes

In this section, we determine the experimental order of convergence (EOC) of specific IMEX upwind gSBP schemes, particularly for coarse time grids. Regarding the designed spatial and temporal order of the schemes, we study the cases $N \leqslant 3$ and utilize the IMEX schemes of the second order (41) and the third order (42).

First, we reconsider the decreasing exact solution (50) of the linear advection-diffusion equation (1). Table 3 lists the $L^2$-errors and the derived EOC for the various schemes applied to this problem. Here, we compute the numeral solution for two different parameter sets of $a = c = 0.1$ and $a = 1, c = 0.1$ until the final time $T = 10$ and use time steps $\Delta t = \mu \Delta x$ for different constants $\mu$ as indicated in the table. In most cases, the predicted convergence rate is experimentally confirmed as the minimum of the expected orders in space and time. One exception is the poor performance in the case of $\left(\theta_{\mathtt{adv}}, \theta_{\mathtt{diff}}\right) = \left(\frac{1}{2}, 0\right)$, i.e., upwind advection fluxes paired with the BR1 diffusion scheme as this choice suffers from the predicted loss of enhanced IMEX stability for larger time steps, particularly in the third-order case. A second exception is the case of $\left(\theta_{\mathtt{adv}}, \theta_{\mathtt{diff}}\right) = (0, 0)$, which yields a stable scheme but is subject to order reduction regarding the formally second-order scheme. As a matter of fact, this type of order reduction is a known phenomenon which is prominent in DG schemes employing numerical fluxes of central type in the case of odd polynomial degrees, also observed for the central-type kinetic energy preserving numerical fluxes for the compressible fluid flow, e.g., in [16, 26].

**Table 3** Accuracy comparison of IMEX upwind gSBP schemes for the exponentially decreasing solution (50). Computations carried out until the final time $T = 10$

| $(a,c)$, $\Delta t$ | $K$ | Upwind gSBP discretization $(\theta_{\mathtt{adv}}, \theta_{\mathtt{diff}})$ | | | | | | | |
|---|---|---|---|---|---|---|---|---|---|
| | | $\left(\frac{1}{2},\frac{1}{2}\right)$ | | $\left(\frac{1}{2},0\right)$ | | $\left(\frac{1}{4},\frac{1}{4}\right)$ | | $(0,0)$ | |
| | | $L^2$-error | EOC | $L^2$-error | EOC | $L^2$-error | EOC | $L^2$-error | EOC |
| Second-order IMEX upwind gSBP schemes | | | | | | | | | |
| $a = 0.1$ | 20 | 1.01E−01 | – | 1.73E−01 | – | 1.01E−01 | – | 1.46E−01 | – |
| $c = 0.1$ | 40 | 2.31E−02 | 2.13 | – | – | 2.37E−02 | 2.09 | 5.07E−02 | 1.53 |
| $N = 1$ | 80 | 6.42E−03 | 1.85 | – | – | 6.60E−03 | 1.84 | 2.27E−02 | 1.16 |
| $\Delta t = 25\Delta x$ | 160 | 1.55E−03 | 2.05 | – | – | 1.60E−03 | 2.04 | 1.09E−02 | 1.06 |
| | 320 | 3.81E−04 | 2.02 | – | – | 3.96E−04 | 2.01 | 5.40E−03 | 1.01 |
| $a = 0.1$ | 20 | 2.10E−02 | – | 1.68E−02 | – | 2.27E−02 | – | 8.62E−02 | – |
| $c = 0.1$ | 40 | 5.72E−03 | 1.88 | 4.41E−03 | 1.93 | 7.28E−03 | 1.64 | 4.31E−02 | 1.00 |
| $N = 1$ | 80 | 1.49E−03 | 1.94 | 1.12E−03 | 1.98 | 2.06E−03 | 1.82 | 2.15E−02 | 1.00 |
| $\Delta t = 5\Delta x$ | 160 | 3.82E−04 | 1.96 | 2.81E−04 | 1.99 | 5.49E−04 | 1.91 | 1.08e−02 | 0.99 |
| | 320 | 9.68E−05 | 1.98 | 7.04E−05 | 2.00 | 1.42E−04 | 1.95 | 5.39E−03 | 1.00 |
| $a = 1$ | 20 | 7.93E−02 | – | 7.66E−02 | – | 7.21E−02 | – | 5.61E−02 | – |
| $c = 0.1$ | 40 | 2.03E−02 | 1.97 | 2.01E−02 | 1.93 | 1.98E−02 | 1.86 | 2.79E−02 | 1.01 |
| $N = 1$ | 80 | 5.12E−03 | 1.99 | 5.10E−03 | 1.98 | 5.06E−03 | 1.97 | 1.39E−02 | 1.01 |
| $\Delta t = 0.5\Delta x$ | 160 | 1.29E−03 | 1.99 | 1.28E−03 | 1.99 | 1.29E−03 | 1.97 | 6.97E−03 | 1.00 |
| | 320 | 3.25E−04 | 1.99 | 3.20E−04 | 2.00 | 3.31E−04 | 1.96 | 3.48E−03 | 1.00 |
| Third-order IMEX upwind gSBP schemes | | | | | | | | | |
| $a = 0.1$ | 20 | 4.40E−04 | – | 3.70E−04 | – | 3.93E−04 | – | 3.80E−04 | – |
| $c = 0.1$ | 40 | 6.12E−05 | 2.85 | 5.06E−05 | 2.87 | 5.36E−05 | 2.87 | 5.04E−05 | 2.91 |
| $N = 2$ | 80 | 7.90e−06 | 2.95 | 6.54E−06 | 2.95 | 6.92E−06 | 2.95 | 6.56E−06 | 2.94 |
| $\Delta t = 5\Delta x$ | 160 | 1.02E−06 | 2.95 | 1.17E−05 | – | 8.90E−07 | 2.96 | 8.40E−07 | 2.97 |
| | 320 | 1.28E−07 | 2.99 | – | – | 1.12E−07 | 2.99 | 1.06E−07 | 2.99 |
| $a = 1$ | 20 | 6.13E−04 | – | – | – | 6.33E−04 | – | 6.74E−04 | – |
| $c = 0.1$ | 40 | 7.94E−05 | 2.95 | – | – | 8.08E−05 | 2.97 | 8.34E−05 | 3.01 |
| $N = 2$ | 80 | 1.05E−05 | 2.92 | – | – | 1.04E−05 | 2.96 | 1.05E−05 | 2.99 |
| $\Delta t = 0.5\Delta x$ | 160 | 1.37E−06 | 2.94 | – | – | 1.33E−06 | 2.97 | 1.32E−06 | 2.99 |
| | 320 | 1.76E−07 | 2.96 | – | – | 1.68E−07 | 2.98 | 1.65E−07 | 3.00 |

Furthermore, as in [38], we now supplement the linear advection-diffusion equation (1) by a source term, i.e., we consider the modified advection-diffusion problem

$$\begin{cases} u_t + u_x &= cu_{xx} + g(x,t), \qquad (x,t) \in (-\pi,\pi) \times (0,T), \\ u(x,0) &= \sin x, x \in (-\pi,\pi). \end{cases} \tag{51}$$

The source term is set to $g(x,t) = \mathrm{e}^{ct}(2c\sin x + \cos x)$ to provide the exact solution $u(x,t) = \mathrm{e}^{ct}\sin x$, which is exponentially growing in time.

Computations are carried out for a diffusion parameter of $c = 0.1$ until the final time $T = 10$ with time steps $\Delta t = \mu\Delta x$ with $\mu$ as indicated in Table 4. Again, the anticipated order of convergence is reached in all setups with two exceptions. Analogously to the exponentially decreasing solution, order reduction occurs for the central scheme $(\theta_{\mathtt{adv}}, \theta_{\mathtt{diff}}) = (0,0)$ for $N = 1$. Furthermore, the pairing of upwind advection with BR1 diffusion replicates the loss of the enhanced IMEX stability for all but the first time step choice in Table 4.

## 4.3 Numerical Results for the Viscous Burgers' Equation

Finally, we study the behavior of IMEX upwind gSBP schemes with respect to the viscous Burgers' equation

$$u_t(x,t) + \left(\frac{1}{2}u^2(x,t)\right)_x = cu_{xx}(x,t),\ c = 0.1,$$

as a prototype containing nonlinear convective terms, supplemented by the initial condition

$$u(x,0) = u_0(x) = \sin x, \quad x \in [-\pi,\pi]$$

and periodic boundary conditions. Figure 1 shows the initial solution as well as the numerical solution for $c = 0.1$ at a given output time $T$. Two different upwind gSBP schemes with $N = 2$ are employed within the following semi-discrete formulation of the viscous Burgers' equation,

$$\frac{\mathrm{d}\mathbf{u}}{\mathrm{d}t} + \left(\frac{\mathbf{D}^+ + \mathbf{D}^-}{2}\right)\mathbf{f} - \|\mathbf{u}\|_\infty \mathbf{M}^{-1}\mathbf{C}\mathbf{u} = c\mathbf{D}_2\mathbf{u} \tag{52}$$

with the vector of flux values $\mathbf{f} = \frac{1}{2}\mathbf{u}^2$ and the matrix $\mathbf{C}$ as in (12). This semi-discrete form is similar to the one used in [34] for the inviscid Burgers' equation but is based on the equation in a conservative form instead of the skew-symmetric formulation employed in [34]. If $\mathbf{D}^-$ is constructed from a DG scheme using upwind advection fluxes, the semi-discretization (52) is equivalent to the use of global Lax-Friedrichs (LF) fluxes for the nonlinear convective term and will be labeled accordingly. If $\mathbf{D}^-$ is constructed from central fluxes, we have $\mathbf{C} = \mathbf{0}$ and thus obtain a central discretization of the inviscid fluxes. For time integration, the second-order IMEX scheme (41) is applied with a time step size of $\Delta t = 0.1$. This same time step is used on two different grids with $K = 50$ and $K = 100$ elements, respectively. First, we compute approximate solutions with the upwind gSBP scheme consisting of BR1 diffusion paired with global LF fluxes, i.e., the first of the above mentioned variants for $\mathbf{D}^-$. Since for linear advection, global LF fluxes reduce to upwind fluxes, the global LF fluxes can be viewed as a generalization of upwind fluxes to the nonlinear case. For this setup, the simulations needs to be stopped at $T = 1.8$ and $T = 1$ due to the expected lack of the stability, as shown in the top row of Fig. 1. These results indicate that also for nonlinear problems, a pairing of the BR1 scheme with upwind-type fluxes might be less favorable in terms of the stability in the context of IMEX advection-diffusion splitting. A remedy is provided in the bottom row of Fig. 1 using central fluxes together with BR1. The corresponding stable numerical solution is shown at output time $T = 2$. This solution is visually indistinguishable from a reference solution computed on a fine grid of $K = 1\,000$ elements using the upwind gSBP scheme defined by $(\theta_{\mathtt{adv}}, \theta_{\mathtt{diff}}) = \left(\frac{1}{2}, \frac{1}{2}\right)$ and $N = 3$ combined with third-order IMEX time integration using a small time step of $\Delta t = 0.01$. In addition, refining the grid does not necessitate choosing smaller time steps for this combination of central first-derivative and second-derivative gSBP operators.

**Table 4** Accuracy comparison of IMEX upwind gSBP schemes for the exponentially growing solution (51) with $c = 0.1$. Computations carried out until the final time $T = 10$

| $N, \Delta t$ | $K$ | Upwind gSBP discretization $(\theta_{\mathtt{adv}}, \theta_{\mathtt{diff}})$ | | | | | | | |
|---|---|---|---|---|---|---|---|---|---|
| | | $\left(\frac{1}{2}, \frac{1}{2}\right)$ | | $\left(\frac{1}{2}, 0\right)$ | | $\left(\frac{1}{4}, \frac{1}{4}\right)$ | | $(0, 0)$ | |
| | | $L^2$-error | EOC | $L^2$-error | EOC | $L^2$-error | EOC | $L^2$-error | EOC |
| Second-order IMEX time integration | | | | | | | | | |
| $N = 1$ | 20 | 9.24E−02 | – | 1.62E−01 | – | 1.44E−01 | – | 8.25E−01 | – |
| $\Delta t = \Delta x$ | 40 | 2.24E−02 | 2.04 | 4.01E−02 | 2.01 | 2.54E−02 | 2.50 | 4.14E−01 | 0.99 |
| | 80 | 6.42E−03 | 1.80 | 1.01E−02 | 1.99 | 5.79E−03 | 2.13 | 2.07E−01 | 1.00 |
| | 160 | 1.87E−03 | 1.78 | 2.53E−03 | 2.00 | 2.13E−03 | 1.44 | 1.04E−01 | 0.99 |
| | 320 | 5.17E−04 | 1.85 | 6.37E−04 | 1.99 | 7.14E−04 | 1.58 | 5.18E−02 | 1.01 |
| $N = 2$ | 20 | 6.72E−04 | – | – | – | 1.10E−03 | – | 1.87E−03 | – |
| $\Delta t = 0.5\Delta x$ | 40 | 1.61E−04 | 2.06 | – | – | 1.90E−04 | 2.53 | 2.39E−04 | 2.97 |
| | 80 | 3.42E−05 | 2.23 | – | – | 3.24E−05 | 2.55 | 3.35E−05 | 2.83 |
| | 160 | 6.37E−06 | 2.42 | – | – | 5.83E−06 | 2.47 | 5.70E−06 | 2.56 |
| | 320 | 1.28E−06 | 2.32 | – | – | 1.22E−06 | 2.26 | 1.20E−06 | 2.25 |
| $N = 3$ | 20 | 1.04E−04 | – | – | – | 1.03E−04 | – | 2.54E−04 | – |
| $\Delta t = 0.3\Delta x$ | 40 | 2.57E−05 | 2.02 | – | – | 2.57E−05 | 2.00 | 3.87E−05 | 2.71 |
| | 80 | 6.43E−06 | 2.00 | – | – | 6.43E−06 | 2.00 | 7.37E−06 | 2.39 |
| | 160 | 1.61E−06 | 2.00 | – | – | 1.61E−06 | 2.00 | 1.67E−06 | 2.14 |
| | 320 | 4.02E−07 | 2.00 | – | – | 4.02E−07 | 2.00 | 4.06E−07 | 2.04 |
| Third-order IMEX time integration | | | | | | | | | |
| $N = 2$ | 20 | 5.55E−04 | – | – | – | 1.04E−03 | – | 1.83E−03 | – |
| $\Delta t = 0.5\Delta x$ | 40 | 1.41E−04 | 1.98 | – | – | 1.74E−04 | 2.58 | 2.26E−04 | 3.02 |
| | 80 | 2.90E−05 | 2.28 | – | – | 2.68E−05 | 2.70 | 2.81E−05 | 3.01 |
| | 160 | 4.52E−06 | 2.68 | – | – | 3.72E−06 | 2.85 | 3.51E−06 | 3.00 |
| | 320 | 6.27E−07 | 2.85 | – | – | 4.89E−07 | 2.93 | 4.39E−07 | 3.00 |
| $N = 3$ | 20 | 1.74E−05 | – | – | – | 1.46E−05 | – | 2.32E−04 | – |
| $\Delta t = 0.3\Delta x$ | 40 | 1.82E−06 | 3.26 | – | – | 1.63E−06 | 3.16 | 2.90E−05 | 3.00 |
| | 80 | 1.65E−07 | 3.46 | – | – | 2.00E−07 | 3.03 | 3.62E−06 | 3.00 |
| | 160 | 1.47E−08 | 3.49 | – | – | 1.84E−08 | 3.44 | 4.53E−07 | 3.00 |
| | 320 | 1.49E−09 | 3.30 | – | – | 1.71E−09 | 3.43 | 5.66E−08 | 3.00 |

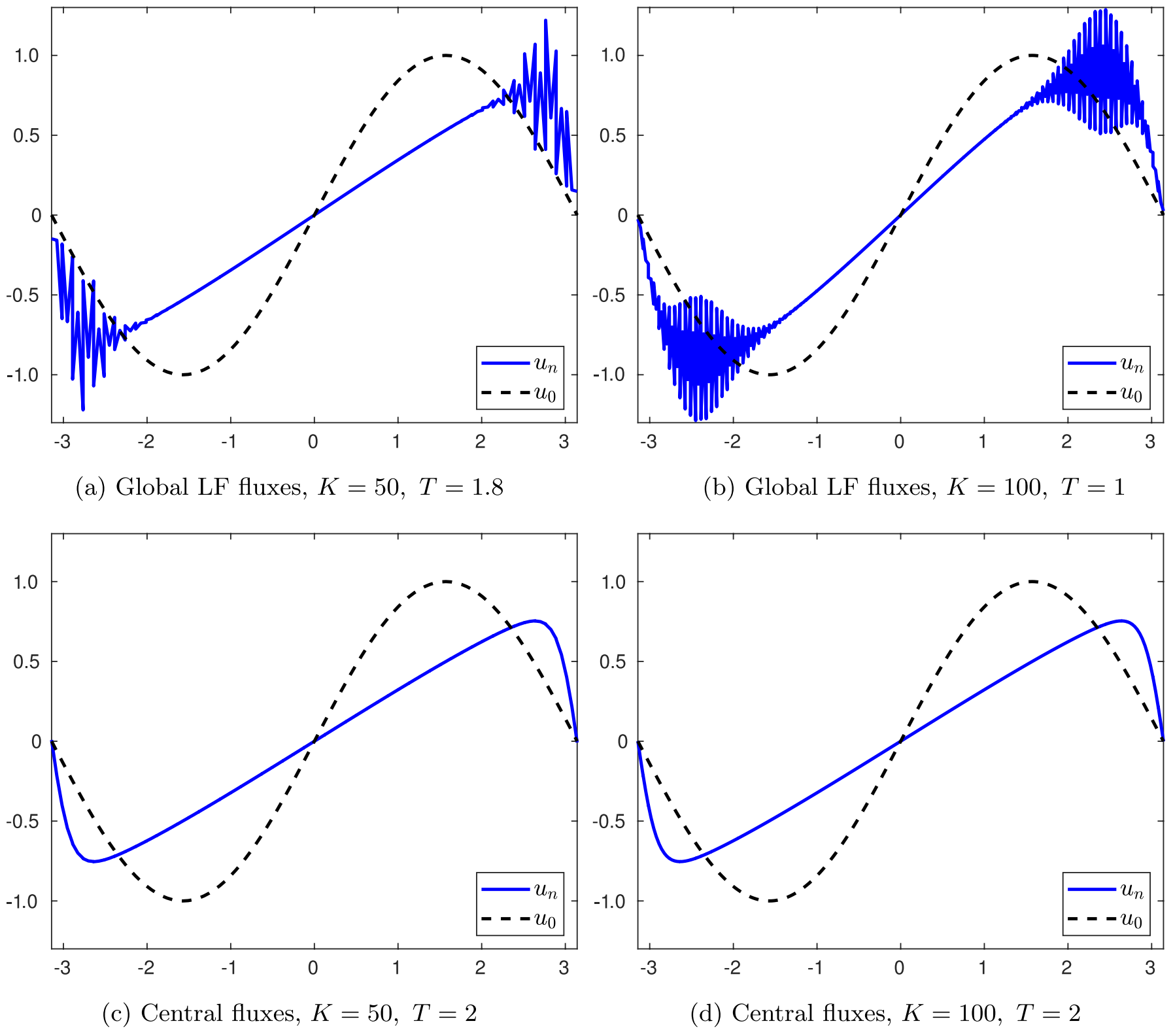

**Fig. 1** Third-order IMEX upwind gSBP approximation of the viscous Burgers' equation until the final time $T = 2$ using BR1 diffusion fluxes ($\theta_{\mathtt{diff}} = 0$) and upwind fluxes for advection ($\theta_{\mathtt{adv}} = \frac{1}{2}$, top row) compared to central fluxes ($\theta_{\mathtt{adv}} = 0$, bottom row)

## 5 Conclusion and Outlook

The current fully discrete energy stability analysis for linear advection-diffusion problems discretized by upwind gSBP schemes in space and IMEX-RK schemes in time based on advection-diffusion IMEX splitting demonstrates that if the employed upwind gSBP operators are chosen to fulfill $\mathbf{D}_2 = \mathbf{D}^{-}\mathbf{D}^{+}$, the allowable time step size is independent of grid refinement, although the advective terms are discretized explicitly. In one space dimension, upwind gSBP schemes represent a general framework including standard DG schemes, thus the present results also hold for the LDG scheme combined with upwind advection fluxes. Furthermore, the analysis in this work is based on the discrete energy provided by the corresponding SBP norm matrix and does not involve inverse constants depending on the discretization order in space; therefore, the resulting time step bound is independent of the polynomial degree of the approximate solution. While previous work for DG schemes has demonstrated that the combination of upwind advection fluxes with the BR1 diffusion scheme does not provide the unconditional stability in the sense of grid-independent time step restrictions, the current work shows that central advection fluxes are compatible with BR1. However, for this kind of overall central scheme, additional artificial dissipation might be necessary for discontinuous initial conditions and in the nonlinear convection-dominated case. These additional artificial dissipation terms could again be discretized either implicitly or explicitly within an IMEX approach. The numerical results for the nonlinear viscous Burgers' equation in this work indicate that a compatible pairing of first-derivative and second-derivative upwind gSBP operators is also required for nonlinear conservation laws with dissipative terms, if the enhanced IMEX stability is desired. This should be investigated more closely with respect to a discrete energy analysis which might require skew-symmetric formulations as in [13, 34] instead of the divergence form of the given conservation law.

Although the focus of this work is restricted to one-dimensional advection-diffusion problems, the presented theoretical results on the enhanced IMEX stability can be carried over to spatial discretizations in higher dimensions as long as the upwind gSBP property is fulfilled and the corresponding first-derivative and second-derivative operators are compatible. This extension relies on the availability of such upwind gSBP operators on structured or unstructured meshes. While structured meshes can be dealt with using tensor-product grids, multi-dimensional generalized SBP operators on simplex meshes have been constructed in [18]. These gSBP operators fulfill a multi-dimensional analog of the properties in Definition 1 with $\mathbf{C} = \mathbf{0}$, i.e., without upwind dissipation. Considering for instance a two-dimensional linear advection-diffusion problem of the form

$$u_t + u_x + u_y = c\,\Delta u, \; c > 0$$

subject to periodic boundary conditions, space discretization may be carried out by first-derivative upwind gSBP operators $\mathbf{D}_x^-$ and $\mathbf{D}_y^-$ discretizing $\frac{\partial}{\partial x}$ and $\frac{\partial}{\partial y}$, respectively, and compatible second-order derivative operators $\mathbf{D}_{2,x} = \mathbf{D}_x^-\mathbf{D}_x^+$ and $\mathbf{D}_{2,y} = \mathbf{D}_y^-\mathbf{D}_y^+$. Due to the upwind gSBP property (11), Proposition 2 then extends to the two-dimensional case and the extension of the stability results of Theorems 2 and 3 is only technical, since it involves an additive decomposition into the terms in $x$- and $y$-directions. The construction of different multi-dimensional upwind gSBP operators with $\mathbf{C} \neq \mathbf{0}$ and their application and comparison in the context of practically relevant multi-dimensional advection-diffusion problems seems to be a promising approach for future exploration. This investigation would also profit from a comparison of the upwind gSBP approach to the multi-dimensional IMEX-LDG schemes considered in [40].

**Funding** Open Access funding enabled and organized by Projekt DEAL.

**Data Availability** Data/code associated with this article is available on request from the author.

## Compliance with Ethical Standards

**Conflict of Interest** The author declares that there is no conflict of interest.